\documentclass[preprint,12pt]{elsarticle}
\usepackage[english]{babel}
\usepackage[utf8x]{inputenc}
\usepackage[T1]{fontenc}
\usepackage{amsmath,amssymb,amsfonts}
\usepackage{color}
\usepackage{float}
\usepackage{graphicx}
\usepackage{tablefootnote}
\usepackage{tabu}
\usepackage{footnote}
\usepackage{mwe}
\usepackage{subfigure}
\usepackage{etex}
\usepackage{caption}
\usepackage{booktabs}
\usepackage{multicol}
\usepackage{multirow}
\usepackage{lipsum}
\usepackage{threeparttable}
\usepackage{scalerel,stackengine}
\usepackage{epstopdf}
\usepackage[toc,page]{appendix}
\usepackage{fixmath}
\usepackage{tikz}
\usepackage{makecell}
\usetikzlibrary{intersections,shapes.arrows,calc,automata,shapes}
\usetikzlibrary{arrows,decorations.markings}
\usetikzlibrary{plotmarks}
\usepackage{pgfplots}
\pgfplotsset{compat=1.14}
\usepackage{lineno}
\usepackage{changes}
\usepackage{titlesec}
\usepackage{hyperref}
\usepackage{url}
\usepackage{breqn}
\usepackage{graphicx,color,colordvi}
\usepackage{latexsym,amsbsy,epsfig,fancybox}
\usepackage[margin=2.5cm]{geometry}
\hypersetup{
    colorlinks=true,
    linkcolor=blue,
}
\usepackage{float}
\usepackage{mathrsfs}
\usepackage{mathtools}
\usepackage{stmaryrd,dsfont}
\usepackage{isomath}
\usepackage{marvosym}
\usepackage{gensymb}
\usepackage{bm,xstring}
\usepackage{etex}
\usepackage{enumitem}
\newlist{steps}{enumerate}{1}
\setlist[steps, 1]{label = Step \arabic*:}
\usetikzlibrary{shapes.misc}
\graphicspath{{figures/}}
\tikzset{cross/.style={cross out, draw=black, minimum size=2*(#1-\pgflinewidth), inner sep=0pt, outer sep=0pt},
cross/.default={1pt}}

\newcommand{\eref}[1]{Equation~(\ref{#1})}

\newcommand{\fref}[1]{Figure~\ref{#1}}

\newcommand{\tref}[1]{Table~\ref{#1}}

\newcommand{\markerfour}{\raisebox{0.5pt}{\tikz{\node[draw,scale=0.4,regular polygon, regular polygon sides=4,fill=none](){};}}}

\newcommand{\xx}{\mathbf{x}}

\stackMath
\newcommand\reallywidehat[1]{%
\savestack{\tmpbox}{\stretchto{%
  \scaleto{%
    \scalerel*[\widthof{\ensuremath{#1}}]{\kern-.6pt\bigwedge\kern-.6pt}%
    {\rule[-\textheight/2]{1ex}{\textheight}}
  }{\textheight}%
}{0.5ex}}%
\stackon[1pt]{#1}{\tmpbox}%
}

\tikzstyle{nicebox}=[draw=black!100, fill=white!10, rectangle, inner sep=4pt, inner ysep=16pt]
\tikzstyle{niceboxtitle}=[draw=black!100, fill=white, text=black, rectangle]

\usetikzlibrary{arrows,decorations.markings}
\biboptions{sort&compress}
\begin{document}
\begin{frontmatter}

\title{
Adaptive Importance Sampling based Neural Network framework for Reliability and Sensitivity Prediction for Variable Stiffness Composite Laminates with hybrid uncertainties}

\author[ind]{Tittu Varghese Mathew}
\author[ind]{P Prajith}
\author[label3]{R. O. Ruiz}
\author[label2]{E. Atroshchenko}
\author[ind]{S Natarajan\corref{cor1}}

\address[ind]{Integrated Modelling and Simulation Lab, Department of Mechanical Engineering, Indian Institute of Technology Madras, Chennai-600036, India.}
\address[label3]{Department of Civil Engineering, Universidad de Chile, Av. Blanco Encalada 2002, Santiago,Chile}
\address[label2]{School of Civil and Environmental Engineering, University of New South Wales, Sydney, Australia}

\cortext[cor1]{Corresponding author}
\fntext[fnlab1]{Department of Mechanical Engineering, Indian Institute of Technology Madras, Chennai-600036, India. E-mail: snatarajan@iitm.ac.in; sundararajan.natarajan@gmail.com}

\begin{abstract}

In this work, we propose to leverage the advantages of both the Artificial Neural Network (ANN) based Second Order Reliability Method (SORM) and Importance sampling to yield an Adaptive Importance Sampling based ANN, with specific application towards failure probability and sensitivity estimates of Variable Stiffness Composite Laminate (VSCL) plates, in the presence of multiple independent geometric and material uncertainties. The performance function for the case studies is defined based on the fundamental frequency of the VSCL plate. The accuracy in both the reliability estimates and sensitivity studies using the proposed method were found to be in close agreement with that obtained using the ANN based brute-force MCS method, with a significant computational savings of 95\%. Moreover, the importance of taking into account the randomness in ply thickness for failure probability estimates is also highlighted quantitatively under the sensitivity studies section. 
\end{abstract}


\end{frontmatter}
\section{Introduction} \label{intro}

None of the above mentioned parameters are deterministic in nature, thanks to the inherent statistical nature of material properties of the constituents and the unavoidable fabrication inaccuracies in ply layup and fiber placements.

Composite materials are used extensively both in primary as well as secondary structures of aerospace and mechanical structures. During the course of a product lifetime, such components are exposed to harsh environments, including mechanical vibration. They are found to exhibit superior strength-to-weight and strength-to-stiffness ratios. Since composites are an amalgamation of two or more constituent materials, its effective dynamic properties are governed by its material mechanical properties, its density, the stacking sequence and thickness of each ply, orientation angle and so on. 

All the above mentioned parameters are not exempt of uncertainties, mainly due to the variability of the material constituents and the unavoidable fabrication inaccuracies in ply layout and fiber placements. In this context, the uncertainty associated to model parameters is understood as incomplete information, facilitating their probabilistic description grounded in the information theory. Hence, all such mechanical and geometric uncertainties must be quantified and taken into account for probabilistic prediction of composite structure response. Natural frequencies of components within safety critical assemblies often tend to be so close to each other, that a small uncertainty in the dynamic characteristics of any of its components can have a catastrophic effect on the overall system performance. Therefore, we must be able to quantify the uncertainty (at least in the first dynamic response) of a critical component made of composites, and as a by-product predict its failure reliability based on a set user-defined criteria. It is also critical to perform a detailed sensitivity study of laminate behaviour to the uncertainty in each input, modeled here as a random variable.

Laminated fibre-reinforced composites \cite{altenbach2011mechanics,reddy2004mechanics} are made up of straight and unidirectional fibers that are homogeneously distributed in each lamina. For such materials, from a macroscopic sense, it is safe to consider that its stiffness does not vary at the laminate level. However, from the early nineties, a shift in the general trend in composite manufacturing has been observed, where the stiffness is purposely made to vary within a lamina, which supposedly lead to more efficient designs \cite{gurdal1993plane,hyer1991use2}. There exists several ways to achieve Variable Stiffness Composite Laminates (VSCL). A few among them are by using curvilinear fibers \cite{tatting1998analysis,wu2006variable,lopes2008variable}, by varying the volume fraction of fibers \cite{setoodeh2006design,senocak2007analysis}, or by adding and dropping plies to the laminates \cite{cairns1999design,her2002stress}. There exists a few advantages for going with variable fibre orientation over others, such as, continuous variation of the stiffness with membrane coordinates \cite{ribeiro2014review} and avoiding abrupt changes in the thickness which may lead to stress concentrations. One other reason that complements the increasing availability of curved fibre laminates is the capability of present tow-placement machines to control the fibre tow placement individually.

The study of static and dynamic response of VSCL plates have received considerable attention. The natural modes of VSCL plates with curvilinear fibers using the $p$-version finite element, that follows a third order deformation theory, was examined in \cite{akhavan2011natural}. Geometric non-linear vibration studies of VSCL plates was studied using the first order shear deformation theory (FSDT) in \cite{ribeiro2012non}. The dynamic instability of VSCL plates using the Rayleigh-Ritz method was studied in \cite{loja2017dynamic}. Anand \textit{et.al.,}~\cite{venkatachari2018variable} studied the free vibration characteristics of VSCL plates using a higher order accurate theory. The influence of environmental effects and geometric discontinuities were also considered. All the aforementioned studies assumed that the geometry and the material parameters were deterministic.

Numerical investigation of the stochastic vibration behaviour of Constant Stiffness Composite Laminate (CSCL) have been studied in~\cite{haldar2000probability,sepahvand2016spectral,chen2008spectral}. A mean-centered second-moment method to study the free vibration and reliability study of composite cantilevers was done in \cite{oh1997free}. A stochastic finite element method for the analysis of composite plates was proposed in \cite{haldar2000probability}. A First Order Reliability Method (FORM) to study shear deformable laminated composite plates was used in \cite{gosling2014high}. For more details, interested readers are referred to a comprehensive review on the reliability studies of composite structures and references therein.

Monte Carlo Simulations (MCS) are another brute-force method for reliability analysis, given a large sample size and at the expense of considerable computational time to achieve a reasonable accuracy. \citet{shaker2008stochastic} was one of the first in using FORM and Second Order Reliability Method (SORM) for reliability analysis of composite plates under free vibration using the stochastic finite element method. However, the uncertainties in that study were limited to only laminae stiffness properties, material density and ply orientation angles. The SORM has been proven to be a superior tool over the MCS for stochastic analysis of free vibrating composite plates. But, one of the major drawbacks of the SORM is that it requires the second order derivative of the performance function w.r.t each of the random variables in the system. This limitation could be partially circumvented using a symbolic solver. However, such a circumvention can become obsolete when the derivative of the performance function w.r.t the random variables like the ply thickness are required. \cite{zhang2016effect} showed the importance of considering uncertainty in the ply thickness for reliability analysis of fiber-reinforced composites under multi-axial loads. 

Structural reliability analysis using Artificial Neural Network (ANN) has been in the mainstream research since past two decades. \citet{hurtado2001neural} gave a comprehensive overview on the use of neural networks in stochastic mechanics. ANN has been successfully used in the structural reliability and optimization studies in the past \cite{papadrakakis1998structural,hurtado2013structural,deng2005structural,schueremans2005benefit,nazari2015three,chok2016neural}. The art of combining ANN with reliability methods like FORM, SORM or MCS is termed as "Neural Network based Reliability". A comparison of ANN against other reliability methods was done in \cite{gomes2004comparison}, wherein a significant reduction in computational time was reported. A comparative study between Multi-Layer Perceptron (MLP) and Radial-Basis Function (RBF) based Neural Network architecture in replacing the FEM in the probabilistic analysis of structures was provided in \cite{hurtado2001neural}. A new, higher order neuron was developed for the deep neural network model to improve its prediction performance for the compressive strength of foamed concrete, another form of composite material, in \cite{nguyen2019deep}. ~\cite{lopes2010reliability,elhewy2006reliability} substituted FEM with a trained neural network to perform reliability studies on composite plates under static loading. 

However, to the best of knowledge of the authors, there exists very few literature that deals with the reliability studies of VSCL plates using FORM/SORM/MCS. \citet{sohouli2017design} presented a detailed reliability study on the VSCL plate structures using two different limit state function definition, one based on the tip deflection and the other based on the first-ply failure criterion by means of Tsai-Wu criterion. However, only uncertainties in the material were taken into account and the analysis was being limited to only static concentrated/distributed load cases. In the present work, the advantage of both the neural network and the adaptive importance sampling in better approximating the performance function at the Maximum Probability of failure Point (MPP) are exploited to predict the failure probability of VSCL w.r.t its first mode of free vibration, while taking into account both material and geometric uncertainties. An MLP-based neural network architecture is used to obtain an ANN-derived function between the input random parameters and the output limit state function values. 

The paper is organized as follows. Section \ref{sec2} presents the finite element formulation for the free vibration studies of a VSCL plate were cutouts are modelled using eXtended Finite Element Method (XFEM). Section \ref{ConvFORMSORM} introduces the reader to the basics of both First Order and Second Order based Reliability Methods and highlight certain limitations in its applications pertaining to the chosen case study using VSCL plates. ANN based surrogate modelling is introduced in Section \ref{sec:4} along with its architecture terminology and the procedure to derive first and second order derivatives of ANN approximated functions. After a short exposition on its extension for deriving FORM/SORM estimates using trained metamodel in Section \ref{ANN_VSCL_prob}, Section \ref{ANN_sampling} discusses about the two main methods adopted for sampling design points for training data, namely the brute-force MCS and the importance sampling. Section \ref{ANN_reliability} concludes the theoretical discussions where global reliability based sensitivity analysis procedure is discussed using ANN trained metamodel using two different approaches. Finally, Section \ref{Numerical_examples} puts all the theoretical formulations into practice using the case study of VSCL plates with a set of carefully chosen random input parameters and their distributions taken mostly from literature. The main take-away points from this work are summarised in Section \ref{Conclusion}.

\section{Theoretical formulation of a VSCL plate} \label{sec2}

\begin{figure*}[h]  
\begin{center}
\includegraphics[width=1\textwidth]{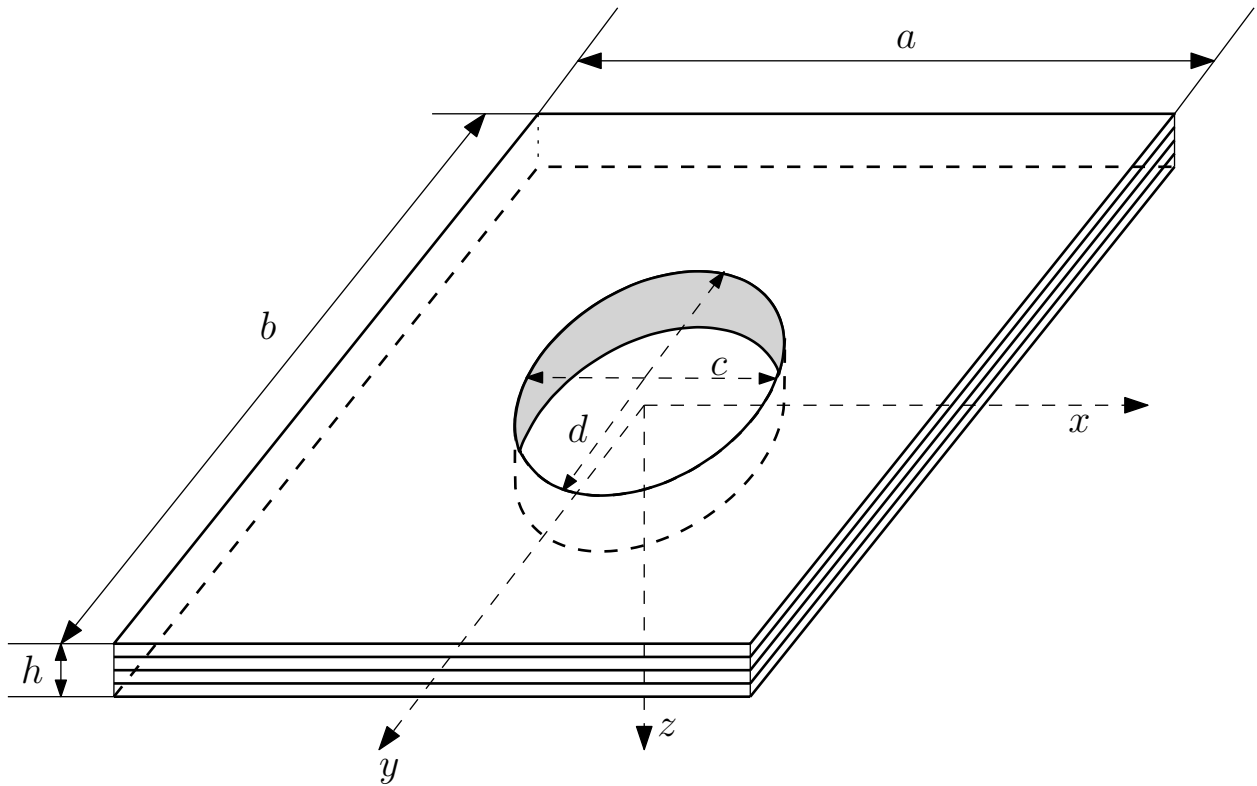}
\end{center}
\captionsetup{justification=centering}
\caption{Composite Laminate plate with an elliptic hole. Cartesian coordinates are shown as $x-y-z$.}
\label{fig:1}
\end{figure*}

The First-order Shear Deformation Theory (FSDT), developed by Reissner and Mindlin for thin and moderately thick plates assumes that there exists a linear variation of the displacement through the thickness of the plate and that the thickness does not change during the deformation. To illustrate the energy formulation for a VSLC plate, a rectangular plate with an elliptic cutout is used, as shown in \fref{fig:1}. 

The plate dimensions are taken as: length $a$ along the $x$ axis, width $b$ along the $y$ axis and thickness $h$ along the $z$ axis. The geometry of the elliptical cutout is represented using the major axis $d$ and the minor axis $c$. The field displacements using the FSDT theory can be expressed as functions of the mid-plane displacements $u_0,v_0,w_0$ and rotations $\beta_x$ and $\beta_y$ of $x-$axis and $y-$axis of the plate respectively as follows: 

\begin{align}
\label{equ:1}
       u(x,y,z) &= u_0(x,y,t) + z \beta_x(x,y,t)  \nonumber \\ 
       v(x,y,z,t) &= v_0(x,y,t) + z \beta_y(x,y,t)  \\ w(x,y,z,t) &= w_0(x,y,t)              \nonumber 
\end{align}

\noindent where $t$ is the time. The linear strain-displacement relationship is given by: 

\begin{equation}
\bm{\varepsilon} = 
\begin{Bmatrix}
\varepsilon_{xx} \\ \varepsilon_{yy}  \\ 2\varepsilon_{xy} \\ 2\varepsilon_{xz} \\ 2\varepsilon_{yz}  
\end{Bmatrix} = 
\begin{Bmatrix}
\frac{\partial u_0}{\partial x} \\
\frac{\partial v_0}{\partial y} \\
\frac{\partial u_0}{\partial y} + \frac{\partial v_0}{\partial x} \\
\frac{\partial w_0}{\partial x} + \beta_x \\
\frac{\partial w_0}{\partial y} + \beta_y  
\end{Bmatrix} + z  
\begin{Bmatrix}
\frac{\partial \beta_x}{\partial x} \\
\frac{\partial \beta_y}{\partial y} \\
\frac{\partial \beta_x}{\partial y} + 
\frac{\partial \beta_y}{\partial x} \\
0 \\ 
0
\end{Bmatrix}
\label{eqn:32}
\end{equation}

\noindent The above strain-displacement relation can be rewritten as follows:  

\begin{equation}
    \bm{\varepsilon} = 
    \begin{Bmatrix}
    \bm{\varepsilon}_p \\ \bm{0}
    \end{Bmatrix} + 
    \begin{Bmatrix}
    z \bm{\varepsilon}_b \\ 
    \bm{\varepsilon}_s
    \end{Bmatrix}
\end{equation}

\noindent where 

\begin{equation}  \label{eqn:strains}
    \bm{\varepsilon}_p = 
    \begin{Bmatrix}
    u_{0,x} \\ v_{0,y} \\ u_{0,y} + v_{0,x}   
    \end{Bmatrix} , 
    \bm{\varepsilon}_b = 
    \begin{Bmatrix}
    \beta_{x,x} \\ \beta_{y,y} \\ \beta_{x,y} + \beta_{y,x} 
    \end{Bmatrix}, 
    \bm{\varepsilon}_s = 
    \begin{Bmatrix}
    \beta_x + w_{0,x} \\ 
    \beta_y + w_{0,y}
    \end{Bmatrix}
\end{equation}

\noindent where $\bm{\varepsilon}_p, \bm{\varepsilon}_b$ and $\bm{\varepsilon}_s$ are the mid-plain strain, bending strain and shear strain. The subscript 'comma' represents the partial derivative w.r.t the spatial coordinate that succeeds it. The membrane stress resultants \textbf{N} and the bending stress resultants \textbf{M} can be related to the membrane strain $\bm\varepsilon_p$ and bending strain $\bm\varepsilon_b$ through the following relation \cite{NATARAJAN20113082}: 

\begin{align}
    &\textbf{N} = 
    \begin{Bmatrix}
    N_{xx} \\ N_{yy} \\ N_{xy} 
    \end{Bmatrix} = \textbf{A} \bm{\varepsilon}_p  +  \textbf{B} \bm{\varepsilon}_b  \nonumber \\ 
    &\textbf{M} = 
    \begin{Bmatrix}
    M_{xx} \\ M_{yy} \\ M_{xy}
    \end{Bmatrix} 
    = \textbf{B} \bm{\varepsilon}_p  +  \textbf{D}_b \bm{\varepsilon}_b 
\end{align}

\noindent where the matrices \textbf{A}$=A_{ij}$, \textbf{B}$=B_{ij}$ and \textbf{D}$_{b}=D_{ij}$ (i,j=1,2,6) are the extensional, bending-extensional and bending stiffness coefficients defined as: 

\begin{equation} \label{eq:ABDmatrix}
    \{A_{ij},B_{ij},D_{ij}\} = \int_{-h/2}^{h/2} \overline{Q}_{ij}\{1,z,z^2 \} dz
\end{equation}

\begin{figure}[htbp!]  
\centering
\includegraphics[width=0.7\textwidth]{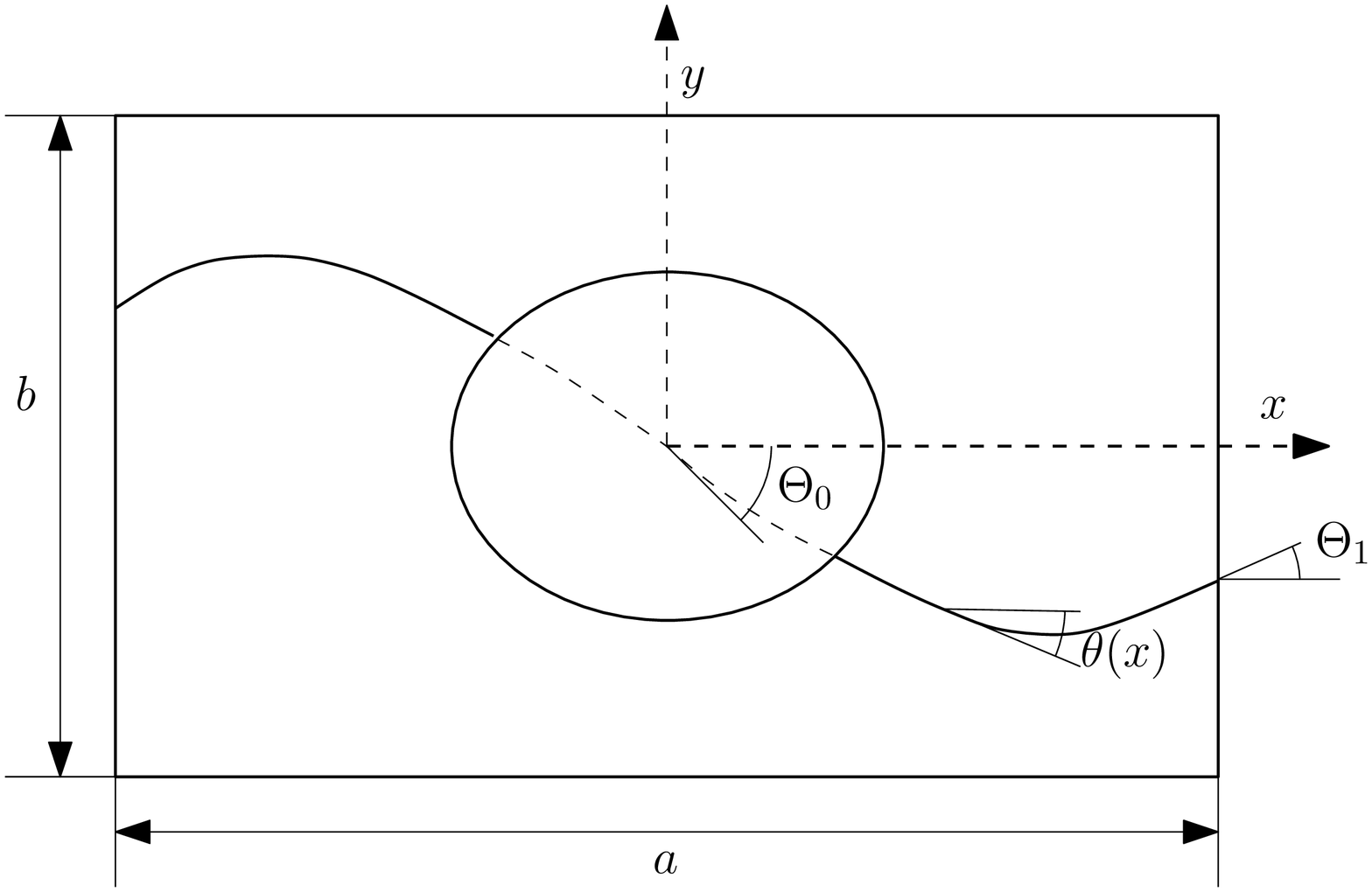}
\captionsetup{justification=centering}
\caption{Curved fiber along with plate dimensions. As shown here, angle $\Theta_0$ has a negative value. Positive value for $\Theta_x$ are defined in the counter-clockwise direction \cite{ribeiro2014review}}
\label{fig:VSCL_sketch}
\end{figure}

\noindent Similarly, the transverse shear force \textbf{Q}$=\{Q_{xz},Q_{yz}\}$ is related to the transverse shear strain $\bm{\varepsilon}_s$ as follows: 

\begin{align}
    &Q_{xz} = K_s \int_{-h/2}^{h/2} \sigma_{xz} dz = K_s \overline{Q}_{55} (\beta_x + w_{0,x}) \nonumber\\ 
    &Q_{yz} = K_s \int_{-h/2}^{h/2} \sigma_{yz} dz = K_s \overline{Q}_{44} (\beta_y + w_{0,y}) 
\end{align}

\noindent where $K_s$ is the shear correction factor. The global stiffness coefficients $\overline{Q}_{ij}$ are derived from the local stiffness coefficients $Q_{ij}$ defined along the local coordinate system that aligns along with the fiber. The local stiffness matrices $Q_{ij}$ are defined as follows: 

\begin{align}
    &Q_{11}=\frac{E_1}{1-\nu_{12}\nu_{21}},Q_{22}=\frac{E_2}{1-\nu_{12}\nu_{21}}, Q_{12}=Q_{21}=\frac{\nu_{12}E_{2}}{1-\nu_{12}\nu_{21}},Q_{16}=Q_{26}=0 \\
    &Q_{44} = G_{23}, Q_{55} = G_{13}, Q_{66} = G_{12}
\end{align}

\noindent where the subscript '1' refers to the local material direction (fiber direction), '2' refers to that direction perpendicular to the local fiber direction and lying within the plane of ply '12'. With the help of a coordinate transformation matrix, the local stiffness coefficients can be transformed to its global stiffness coefficients. However, the curvilinear fibre laminates differ from straight fibre laminates, in the fact that the orientation of the principal material axis does not remain constant in the lamina domain. This means that the transformed stiffness coefficient matrix $\overline{\textbf{Q}}$ doesn't remain constant and become function of spatial $x$ direction. For example, the first element of matrix, namely $\overline{Q}_{11}$ is: 

\begin{equation}
    \overline{Q}_{11} = Q_{11} c(x)^4 + 2 (Q_{12} + 2 Q_{33})c(x)^2 s(x)^2 + Q_{12} s(x)^4
\end{equation}

\noindent where $c(x) = cos(\theta (x))$ and $s(x) = sin(\theta(x))$ in which $\theta(x)$ is the fibre orientation angle of the VSCL plate, measured from the $x$ axis of the plate to the fibre axis angle at position $x$. In this study, it is assumed that the fibre path variation is represented as follows: 

\begin{equation}
    \theta(x) = \frac{2(\Theta_0 - \Theta_1)}{a} |x| + \Theta_0
\end{equation}

\noindent where $\Theta_0$ is the angle between the fibre and the $x$ axis at $x=0$ and $\Theta_1$ is the fibre angle at the panel ends $(x = \pm a/2)$. The fibre orientation for the $k$'th layer is then represented as $<\Theta_0^k | \Theta_1^k>$ as shown in \fref{fig:VSCL_sketch}. However, due to manufacturing constraints, these sets of fiber ply angles cannot take arbitrary combination of values, as was pointed out by Akhavan\textit{.et.al} \cite{akhavan2011natural}. In his seminal work on natural modes of vibration of composite structures with variable fibers, he proposed a relation (c.f.\eref{eq:akhavan}) wherein for a given set of fiber ply pair $<\Theta_0^k | \Theta_1^k>$ to be feasible from a manufacturing perspective, the local curvature of the fiber, $\kappa(x)$, should not exceed a value of 3.28.  

\begin{equation} \label{eq:akhavan}
    \kappa(x) = 2 \frac{(\Theta_1 - \Theta_0)}{a} cos\left( (\Theta_1 - \Theta_0) \frac{x}{a/2} + \Theta_0  \right) < 3.28 
\end{equation}

\noindent The strain energy function $U$ is then given by: 
\begin{align} \label{eqn:PotentialEnergy}
    U(\bm{\delta}) &= \frac{1}{2} \int_{\Omega} \{ \bm{\varepsilon}_{p}^T \textbf{N} + \bm{\varepsilon}_b^T \textbf{M} + \bm{\varepsilon}_s^T \textbf{Q}  \}d\Omega  \nonumber \\
    &= \frac{1}{2} \int_{\Omega} \{ \bm{\varepsilon}_p^T \textbf{A} \bm{\varepsilon}_{p} +  \bm{\varepsilon}_{p}^T \textbf{B} \bm{\varepsilon}_{b} + 
    \bm{\varepsilon}_b^T \textbf{B} \bm{\varepsilon}_{p} +  \bm{\varepsilon}_{b}^T \textbf{D} \bm{\varepsilon}_{b} + 
    \bm{\varepsilon}_s^T \textbf{E} \bm{\varepsilon}_s \} d\Omega
\end{align}

\noindent where $\bm{\delta} = \{u_0, v_0, w_0, \beta_x,\beta_y \}$ is the vector of degrees of freedom associated to the displacement field in finite element discretization. Following the procedure outlined in \cite{rajasekaran1973incremental}, the strain energy function $U(\bm{\delta})$ can be rewritten as: 

\begin{equation} \label{eq:UE}
    U(\bm{\delta}) = \frac{1}{2} \bm{\delta}^T \textbf{K} \bm{\delta}
\end{equation}

\noindent where \textbf{K} is the linear stiffness matrix. In a similar way, the kinetic energy of the plate is given by: 

\begin{equation} \label{eq:KE}
    T(\bm{\delta}) = \frac{1}{2} \int_{\Omega}\{ I_0 (\dot{u}_0^2 + \dot{v}_0^2 + \dot{w}_0^2)+ I_1 (\dot{\theta}_x^2 + \dot{\theta}_y^2) \} d\Omega
\end{equation}

\noindent where $I_0 = \int_{-h/2}^{h/2} \rho dz$ and $I_1 = \int_{-h/2}^{h/2} z^2 \rho dz$, $\rho$ is the mass density. Substituting \eref{eq:UE} and \eref{eq:KE} in a Lagrangian equation of motion, one obtains the following governing equation: 

\begin{equation}
    \textbf{M} \ddot{\bm{\delta}} + \textbf{K} \bm{\delta} = 0
\end{equation}

\noindent where \textbf{M} is the consistent mass matrix. After substituting the characteristic of the time function $\ddot{\bm{\delta}} = \omega^2 \bm{\delta}$, following relation is obtained: 

\begin{equation}
    (\textbf{K} - \omega^2 \textbf{M}) \bm{\delta}= 0 
\end{equation}

\noindent where $\omega$ is the natural frequency. 

\subsection{Spatial discretization}
The plate element employed here is a $\mathcal{C}^o$ continuous shear flexible field consistent element with five degrees of freedom $(u_o,v_o,w_o,\beta_x,\beta_y)$ at four nodes in a 4-noded quadrilateral (QUAD-4) element. The displacement field within the element is approximated by:

\begin{equation}
\{ u_o^e,v_o^e,w_o^e,\beta_x^e,\beta_y^e\} = \sum\limits_{J=1}^4 N_J \{u_{oJ}, v_{oJ}, w_{oJ},\beta_{xJ},\beta_{yJ} \},
\end{equation} 

\noindent where $u_{oJ}, v_{oJ}, w_{oJ},\beta_{xJ},\beta_{yJ}$ are the nodal variables and $N_J$ are the shape functions for the bi-linear QUAD-4 element. In this study, to alleviate the shear locking phenomenon, field redistributed shape functions~\cite{Somashekar1987,natarajanbaiz2011a} are employed.


The finite element framework requires the underlying finite element mesh to conform to the discontinuity surface. The recent introduction of implicit boundary definition-based methods, viz., the extended/generalized FEM (XFEM/GFEM), alleviates the shortcomings associated with the meshing of the discontinuity surface. In this study, the partition of unity framework is employed to represent the discontinuity surface (cutout geometry) independent of the underlying mesh.

\begin{equation}
(u_0^h,v_0^h,w_0^h,\beta_x^h,\beta_y^h)\left(\xx\right) = \underbrace{ \sum_{I \in \mathcal{N}^{\rm{fem}}} N_I(\xx) {(u_0^s,v_0^s,w_0^s,\beta_{x}^s,\beta_{y}^s)}_I}_{\rm FEM} + \underbrace{ \sum_{J \in
\mathcal{N}^{\rm{c}}} N_J(\xx) H(\xx) {(b^{u_0},b^{v_0},b^{w_0}, b^{\beta_x},b^{\beta_y})}_J}_{\rm Enriched~part} 
\label{eqn:platexfem1}
\end{equation}

\noindent where $\mathcal{N}^{\rm{fem}}$ is a set of all the nodes in the finite element mesh and $\mathcal{N}^{\rm{c}}$ is a set of nodes that are enriched with the Heaviside function. In \eref{eqn:platexfem1}, $(u_I^s,v_I^s,w_I^s,\beta_{x_I}^s,\beta_{y_I}^s)$ are the nodal unknown vectors associated with the continuous part of the finite element solution, $b_J$ is the nodal enriched degree of freedom vector associated with the Heaviside (discontinuous) function. In this study, a level set approach is followed to model the cutouts. The geometric interface (for example, the boundary of the cutout) is represented by the zero level curve $\phi \equiv \phi(\xx,t) = 0$. The interface is located from the value of the level set information stored at the nodes. The standard FE shape functions can be used to interpolate $\phi$ at any point $\xx$ in the domain as:

\begin{equation}
\phi(\xx) = \sum\limits_I N_I(\xx) \phi_I
\end{equation}

\noindent where the summation is over all the nodes in the connectivity of the elements that contact $\xx$ and $\phi_I$ are the nodal values of the level set function. For circular cutout, the level set function is given by:

\begin{equation}
\phi_I = || \xx_I - \xx_c|| - r_c
\end{equation}

\noindent where $\xx_c$ and $r_c$ are the center and the radius of the cutout. Based on the level set approach, the Heaviside function $H(\xx)$ in \eref{eqn:platexfem1} becomes: 

\begin{equation}
    H(\xx) = 
    \begin{cases}
    0 ,& \text{if } \phi(\xx) < 0\\
    1  ,& \text{if } \phi(\xx) > 0
    \end{cases}
\end{equation}

\subsubsection{Elements of stiffness and mass matrices}
\begin{figure}[htbp!]  
\centering
\includegraphics[width=0.7\textwidth]{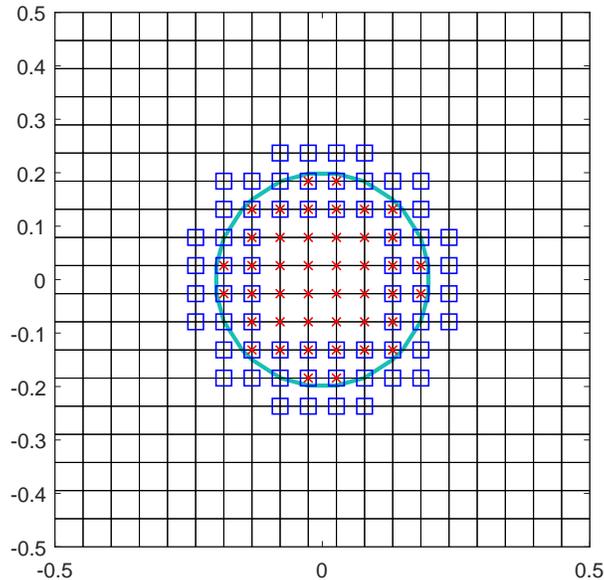}
\captionsetup{justification=centering}
\caption{Mesh for a square plate of dimensions $1 \times 1 $m with a circular cutout of radius $0.2$m, marked with \textcolor{blue}{blue line}. Enriched nodes associated with the elements containing the cutout boundary are shown in blue squares \textcolor{blue}{\protect\markerfour}, while constrained nodes that lie within the cutout are marked as \textcolor{red}{*}.}
\label{fig:XFEM_mesh}
\end{figure}

In this section, the elements of the stiffness and mass matrices are derived using extended finite element formulation. Substituting \eref{eqn:platexfem1} in \eref{eqn:strains}, one gets the following:

\begin{align*}
    &\bm{\varepsilon}_p^h = \sum_{i \in \mathcal{N}^{\rm{fem}}} 
    \begin{bmatrix}
    \frac{\partial N_i}{\partial x} & 0 & 0 & 0 & 0 \\ 
    0 & \frac{\partial N_i}{\partial y} & 0 & 0 & 0 \\ 
    \frac{\partial N_i}{\partial y} & \frac{\partial N_i}{\partial x} & 0 & 0 & 0 
    \end{bmatrix} 
    {\begin{Bmatrix}
    u_0^s \\ v_0^s \\ w_0^s \\ \beta_x^s \\ \beta_y^s
    \end{Bmatrix}}_{i} + 
    \sum_{j \in \mathcal{N}^{\rm{c}}} 
    \begin{bmatrix}
    \reallywidehat{H(\xx)} \frac{\partial N_j}{\partial x} & 0 & 0 & 0 & 0 \\ 
    0 & \reallywidehat{H(\xx)}\frac{\partial N_j}{\partial y} & 0 & 0 & 0 \\ 
    \reallywidehat{H(\xx)}\frac{\partial N_j}{\partial y} & \reallywidehat{H(\xx)}\frac{\partial N_j}{\partial x} & 0 & 0 & 0 
    \end{bmatrix} 
    {\begin{Bmatrix}
    b^{u_0} \\ b^{v_0} \\ b^{w_0} \\ b^{\beta_x} \\ b^{\beta_y}
    \end{Bmatrix}}_{j} \nonumber \\ 
    &\bm{\varepsilon}_b^h = \sum_{i \in \mathcal{N}^{\rm{fem}}} 
    \begin{bmatrix}
    0 & 0 & 0 & \frac{\partial N_i}{\partial x} & 0 \\ 
    0 & 0 & 0 & 0 & \frac{\partial N_i}{\partial y} \\ 
    0 & 0 & 0 & \frac{\partial N_i}{\partial y} & \frac{\partial N_i}{\partial x}
    \end{bmatrix} 
    {\begin{Bmatrix}
    u_0^s \\ v_0^s \\ w_0^s \\ \beta_x^s \\ \beta_y^s
    \end{Bmatrix}}_{i} + 
    \sum_{j \in \mathcal{N}^{\rm{c}}} 
    \begin{bmatrix}
    0 & 0 & 0 & \reallywidehat{H(\xx)}\frac{\partial N_j}{\partial x} & 0 \\ 
    0 & 0 & 0 & 0 & \reallywidehat{H (\xx)}\frac{\partial N_j}{\partial y} \\ 
    0 & 0 & 0 & \reallywidehat{H(\xx)}\frac{\partial N_j}{\partial y} & \reallywidehat{H(\xx)}\frac{\partial N_j}{\partial x}
    \end{bmatrix} 
    {\begin{Bmatrix}
    b^{u_0} \\ b^{v_0} \\ b^{w_0} \\ b^{\beta_x} \\ b^{\beta_y}
    \end{Bmatrix}}_{j} \nonumber \\
    &\bm{\varepsilon}_s^h = \sum_{i \in \mathcal{N}^{\rm{fem}}} 
    \begin{bmatrix}
    0 & 0 & \frac{\partial N_i}{\partial x} & N_i & 0 \\ 
    0 & 0 & \frac{\partial N_i}{\partial y} & 0 & N_i \\ 
    \end{bmatrix} 
    {\begin{Bmatrix}
    u_0^s \\ v_0^s \\ w_0^s \\ \beta_x^s \\ \beta_y^s
    \end{Bmatrix}}_{i} + 
    \sum_{j \in \mathcal{N}^{\rm{c}}} 
    \begin{bmatrix}
    0 & 0 & \reallywidehat{H(\xx)}\frac{\partial N_j}{\partial x} & \reallywidehat{H(\xx)}N_j & 0 \\ 
    0 & 0 & \reallywidehat{H(\xx)}\frac{\partial N_j}{\partial y} & 0 & \reallywidehat{H(\xx)}N_j \\
    \end{bmatrix} 
    {\begin{Bmatrix}
    b^{u_0} \\ b^{v_0} \\ b^{w_0} \\ b^{\beta_x} \\ b^{\beta_y}
    \end{Bmatrix}}_{j} \nonumber \\
\end{align*}
where the derivative of the Heaviside function w.r.t the spatial coordinates exists only at the interface location and $\reallywidehat{H(\xx)} = H(\xx) - H(\mathbf{x}_j)$. The above set of equations can be rewritten in the form of their respective strain displacement matrices as follows: 
\begin{equation} \label{eq:straindispmatrix}
    \bm{\varepsilon}_p^h = \mathbb{B}_p^s ~\bm{\delta}^s + \mathbb{B}_p^e ~\bm{\delta}^e, ~\quad \bm{\varepsilon}_b^h = \mathbb{B}_b^s ~\bm{\delta}^s + \mathbb{B}_b^e ~\bm{\delta}^e, ~\quad 
    \bm{\varepsilon}_s^h = \mathbb{B}_s^s ~\bm{\delta}^s + \mathbb{B}_s^e ~\bm{\delta}^e
\end{equation}
where $\mathbb{B}_p, \mathbb{B}_b$ and $\mathbb{B}_s$ are the membrance, bending and shear strain-displacement matrices, respectively, the superscripts 's' and 'e' stands for the standard and enriched parts. Substituting \eref{eq:straindispmatrix} into \eref{eqn:PotentialEnergy} and upon further simplification, one gets the following elemental stiffness matrix: 
\begin{equation}
    \textbf{K}^e = 
    \begin{bmatrix}
    \textbf{K}^e_{ss} & \textbf{K}^e_{se} \\
    \textbf{K}^e_{es} & \textbf{K}^e_{ee}
    \end{bmatrix}  
\end{equation}
where 
\begin{align*}
    &\textbf{K}_{ss}^e = \int_{\Omega_e} \{{\mathbb{B}_p^s}^T \textbf{A} \mathbb{B}_p^s +  
    {\mathbb{B}_p^s}^T \textbf{B} \mathbb{B}_b^s + 
    {\mathbb{B}_b^s}^T \textbf{B} \mathbb{B}_p^s + 
    {\mathbb{B}_b^s}^T \textbf{D} \mathbb{B}_b^s + 
    {\mathbb{B}_s^s}^T \textbf{E} \mathbb{B}_s^s \} d\Omega \\ 
    &\textbf{K}_{se}^e = \int_{\Omega_e} \{{\mathbb{B}_p^s}^T \textbf{A} \mathbb{B}_p^e +  
    {\mathbb{B}_p^s}^T \textbf{B} \mathbb{B}_b^e + 
    {\mathbb{B}_b^s}^T \textbf{B} \mathbb{B}_p^e + 
    {\mathbb{B}_b^s}^T \textbf{D} \mathbb{B}_b^e + 
    {\mathbb{B}_s^s}^T \textbf{E} \mathbb{B}_s^e \}d\Omega \\
    &\textbf{K}_{es}^e = \int_{\Omega_e} \{{\mathbb{B}_p^e}^T \textbf{A} \mathbb{B}_p^s +  
    {\mathbb{B}_p^e}^T \textbf{B} \mathbb{B}_b^s + 
    {\mathbb{B}_b^e}^T \textbf{B} \mathbb{B}_p^s + 
    {\mathbb{B}_b^e}^T \textbf{D} \mathbb{B}_b^s + 
    {\mathbb{B}_s^e}^T \textbf{E} \mathbb{B}_s^s \}d\Omega \\
    &\textbf{K}_{ee}^e = \int_{\Omega_e} \{{\mathbb{B}_p^e}^T \textbf{A} \mathbb{B}_p^e +  
    {\mathbb{B}_p^e}^T \textbf{B} \mathbb{B}_b^e + 
    {\mathbb{B}_b^e}^T \textbf{B} \mathbb{B}_p^e + 
    {\mathbb{B}_b^e}^T \textbf{D} \mathbb{B}_b^e + 
    {\mathbb{B}_s^e}^T \textbf{E} \mathbb{B}_s^e \}d\Omega 
\end{align*}
Similarly, \eref{eq:KE} can be written as: 
\begin{equation}
    T^e = \frac{1}{2} \{ {\dot{\bm{\delta}}_s}^T, {\dot{\bm{\delta}}_e}^T  \} \textbf{M}^e 
    \begin{Bmatrix}
    \dot{\bm{\delta}}_s \\ \dot{\bm{\delta}}_e 
    \end{Bmatrix}
\end{equation}
Substituting \eref{eqn:platexfem1} into the above relation, one gets to the element mass matrix as follows: 
\begin{equation}
    \textbf{M}^e = 
    \begin{bmatrix}
    \textbf{M}^e_{ss} & \textbf{M}^e_{se} \\
    \textbf{M}^e_{es} & \textbf{M}^e_{ee}
    \end{bmatrix}  = 
    \begin{bmatrix}
    \int_{\Omega_e} \textbf{N}^T \bm{\rho} \textbf{N} d\Omega &  \int_{\Omega_e} \textbf{N}^T \bm{\rho} (H(\xx)\textbf{N}) d\Omega \\ 
    \int_{\Omega_e} (H(\xx)\textbf{N})^T \bm{\rho} \textbf{N} d\Omega & \int_{\Omega_e} (H(\xx)\textbf{N})^T \bm{\rho} (H(\xx)\textbf{N}) d\Omega
    \end{bmatrix}
\end{equation}
where $\bm{\rho}$ consists of $I_0$ and $I_1$. \fref{fig:XFEM_mesh} shows a finite element mesh for a plate with a hole, where both the enriched nodes in the presence of a cutout and constrained nodes are highlighted.


A consequence of adding custom tailored enrichment functions to the FE approximation basis, which are not necessarily smooth functions is that, special care has to be taken in numerically integrating over the elements that are intersected by the discontinuity surface. The standard Gau\ss~quadrature cannot be applied in elements enriched by discontinuous terms, because Gau\ss~quadrature implicitly assumes a polynomial approximation. In the present study, a triangular quadrature with sub-division is employed along with the integration rules described in \tref{table:subcellgausspt}. For the elements that are not enriched, a standard 2 $\times$ 2 Gaussian quadrature rule is used.

\begin{table}[htpb]
\captionsetup{justification=centering}
\caption{Integration rules for enriched and non-enriched elements in the presence of a cutout.}
\centering
\begin{tabular}{lc}
\hline
Element Type & Gau\ss ~points\\
\hline
Non-enriched element & 4  \\
Split element & 3 per triangle  \\
Split blending element & 4 \\
\hline
\end{tabular}
\label{table:subcellgausspt}
\end{table}

\section{Structural Reliability Analysis}\label{ConvFORMSORM}

The analysis and design of complex engineering systems depend heavily on the predictions from numerical models like finite element analysis, while its accuracy depends on how close the digital representation approximates the real-world system. Structural Reliability Analysis (SRA) aims to provide a rational framework to address uncertainties in structural design so that the later can be more objective and less dependent on unrealistic assumptions. The theory behind SRA is formulated around the concept of probability of failure: 

\begin{equation} \label{eq:Pf}
    P_f = \textit{Prob}[g(\mathbf{x}) \leq 0] = \int_{g(\mathbf{x}) \leq 0} f_{X}(\mathbf{x}) \rm{d} \mathbf{x}
\end{equation}

\noindent where $f_{X}(\mathbf{x})$ is the joint Probability Density Function (PDF) of random vector $X$ and $g(\textbf{x})$ is the limit state function, where $g(\textbf{x}) \leq 0$ denotes the failure domain and $g(\textbf{x}) > 0$ is the safe domain. Though the definition in \eref{eq:Pf} looks simple, its exact integration through direct integration is often intractable for two reasons: (1) dimension of integral is usually too high in case of multiple input random parameters and (2) the shape of the limit state function $g(\mathbf{x})$ could be complicated with complex topology, specially for low failure probabilities. SRA methods can be broadly classified into three categories: Taylor series-based approaches such as the First Order Reliability analysis Method (FORM) and Second Order Reliability analysis Method (SORM), simulation based methods such as MCS and its variants and surrogate methods such as the Response Surface Method (RSM) and the Kriging meta-model \cite{aldosary2018structural}. In the following subsections, we will briefly talk about the first two categories followed by the rational in choosing the right approach with our problem in hand. 

\subsection{First Order Reliability Method} \label{FORM}
In this section, a basic introduction to the conventional First Order Reliability Method (FORM) along with its formulation is given. This method is restricted to normal random variable space. So, all the input random variables $x_i$ should be transformed to its reduced variables having a standard normal Gaussian distribution. This step is done so as to simplify the integrand $f_{X}(\mathbf{x})$ in \eref{eq:Pf} so that its contours on the random parameter space will be more regular and symmetric. Assuming that all the $n$ input random variables in $X$ space follows a Gaussian distribution $\mathcal{N}(\mu_i,\sigma_i^2), i = 1,2,\cdots,n$, where $\mu_i$ and $\sigma_i^2$ are the mean and variance of the $i$th input random variable, they can be transformed to a standard normal distribution space $\mathcal{Z} \coloneqq \mathcal{N}(0,1)$ using a linear transformation rule as given below: 

\begin{equation} \label{eq:16}
    z_i = \frac{x_i - \mu_i}{\sigma_i}, \qquad i=1,2,\cdots,n
\end{equation}

\noindent After the transformation, the failure probability integration becomes 

\begin{equation}
    P_f = P(g(\mathbf{z})<0) = \int_{g(\mathbf{z})<0} \phi_{\mathcal{Z}}(\mathbf{z}) \rm{d} \mathbf{z}
\end{equation}

\noindent where $\phi_{\mathcal{Z}}(\mathbf{z})$ is the joint pdf of all input variables in standard normal space.  
\noindent Once the transformation is done, the next step is to approximate the integration boundary $g(\mathbf{z})=0$. FORM is based on the first order Taylor series expansion of the integration boundary as shown below:

\begin{equation} \label{eq:19}
    g(\mathbf{z}) \approx g(\mathbf{z}^*)+\nabla g(\mathbf{z}^*)(\mathbf{z} - \mathbf{z}^*)^T
\end{equation}

\noindent where $\mathbf{z}^*$ is the expansion point and $\nabla g(\mathbf{z}^*)$ is the gradient of $g(\mathbf{z})$ at $\mathbf{z}^*$, defined as

\begin{equation}\label{eq:dGdzi}
    \nabla g(\mathbf{z}^*) = \left( \frac{\partial g(\mathbf{z})}{\partial z_1}, \frac{\partial g(\mathbf{z})}{\partial z_2},\cdots,\frac{\partial g(\mathbf{z})}{\partial z_n}   \right)\bigg\rvert_{\mathbf{z}^*}
\end{equation}

\noindent The point on $g(\mathbf{z})=0$ that has the highest joint probability density is termed as the Most Probable failure Point (MPP). In other words, maximizing the joint pdf $\phi_{\mathcal{Z}}(\mathbf{z})$ along the limit state function $g(\mathbf{z})=0$ gives the location of MPP. The above statement can be mathematically formulated as follows: 

\begin{equation} \label{eq:21}
\max_{\mathbf{z}}~\prod_{i=1}^{n} \frac{1}{\sqrt{2\pi}} \text{exp}\left(- \frac{1}{2}z_i^2  \right), ~~ \textrm{s.t.}~~ g(\mathbf{z}) = 0
\end{equation}

\begin{figure}[htbp!]  
\centering
\includegraphics[width=0.5\textwidth]{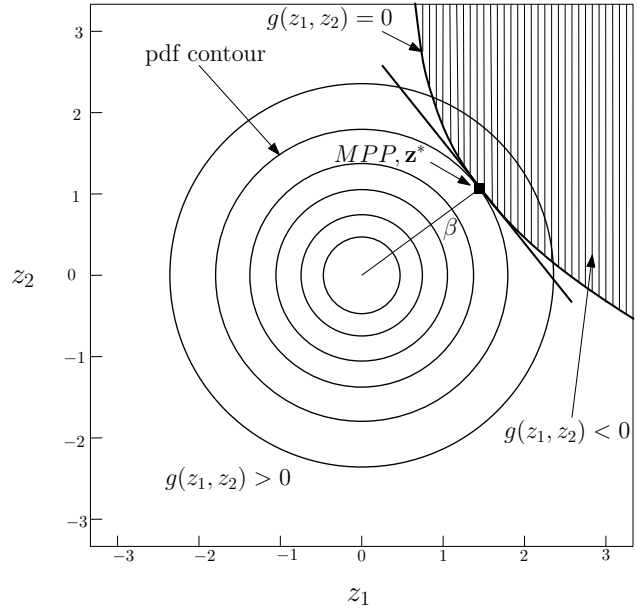}
\captionsetup{justification=centering}
\caption{Probability integration using FORM}
\label{fig:Prob_intg_FORM}
\end{figure}

Since the function to be maximized is equivalent to minimizing $\sum_{i=1}^{n}z_i^2$, the model for the MPP (\eref{eq:21}) can be re-written as follows: 

\begin{equation} \label{eq:22}
\min_{\mathbf{z}} \left\lVert \mathbf{z} \right\rVert \\, ~~ \textrm{s.t.} ~~ g(\mathbf{z}) = 0
\end{equation}

\noindent where $\left\lVert \cdot \right\rVert$ stands for the norm of a vector. The solution to the optimization problem given in \eref{eq:22} is the MPP and its denoted by $\mathbf{z}^* = (z_1^*,\cdots,z_n^*)$. As one can see from \fref{fig:Prob_intg_FORM} for a bivariate case, MPP is the point along $g(z_1,z_2)=0$ which has the shortest distance to the origin of the $\mathcal{Z}-$ space. Once the MPP is computed, the reliability index can be computed as $\beta = \left\lVert \mathbf{z}^*\right\rVert$.

\noindent While the procedure to determine the MPP will be detailed in Section \ref{MPP}, further details in computing the reliability index using FORM is detailed henceforth. Since at MPP $\mathbf{z}^*$, $g(\mathbf{z})=0$, \eref{eq:19} becomes 

\begin{equation} \label{eq:23}
    L(\mathbf{z}) = \sum_{i=1}^{n} \frac{\partial g(\mathbf{z})}{\partial z_i}\bigg\rvert_{\mathbf{z}^*}(z_i - z_i^*) = a_0 + \sum_{i=1}^{n}a_i z_i
\end{equation}

\noindent where 

\begin{equation}
    a_0 = - \sum_{i=1}^{n} \frac{\partial g(\mathbf{z})}{\partial z_i} \bigg\rvert_{\mathbf{z}^*} z_i^*, \quad a_i = \frac{\partial g(\mathbf{z})}{\partial z_i}\bigg\rvert_{\mathbf{z}^*}
\end{equation}

\noindent Since $L(\mathbf{z})$ is a linear function of standard normal variables as per \eref{eq:23}, it is also normally distributed with a mean of $\mu_{L} = a_0$ and standard deviation given by 

\begin{equation}
    \sigma_L = \sqrt{\sum_{i=1}^{n} a_i^2} = \sqrt{\sum_{i=1}^{n} \left( \frac{\partial g}{\partial z_i} \bigg\rvert_{z_i^*}\right)^2}
\end{equation}

\noindent Now using these two information, the probability of failure can be computed as 

\begin{equation}
    P_f \approx P(L(\mathbf{z})<0)) = \Phi \left( \frac{-\mu_L}{\sigma_L} \right)=\Phi \left( \frac{\sum_{i=1}^{n} \frac{\partial g}{\partial z_i}\bigg\rvert_{\mathbf{z}^*} z_i^*}{\sqrt{\sum_{i=1}^{n} \left( \frac{\partial g}{\partial z_i}\bigg\rvert_{\mathbf{z}^*}\right)^2}} \right) = \Phi \left( \sum_{i=1}^{n} \alpha_i z_i^*  \right)
\end{equation}

\noindent where 

\begin{equation}
    \alpha_i = \frac{\frac{\partial g}{\partial z_i}\bigg\rvert_{\mathbf{z}^*}}{\sqrt{\sum_{i=1}^n \left( \frac{\partial g}{\partial z_i} \bigg\rvert_{\mathbf{z}^*}\right)^2}}
\end{equation}

\noindent Let the collection of vectors $\alpha_i$ be represented as follows

\begin{equation}
    \mathbf{a} = (\alpha_1,\alpha_2,\cdots,\alpha_n) = \frac{\nabla g(\mathbf{z}^*)}{\left\lVert \nabla g(\mathbf{z}^*) \right\rVert }
\end{equation}

\noindent The probability of failure can be finally written in the following format: 

\begin{equation} \label{eq:30}
    P_f \approx \Phi \left( \sum_{i=1}^{n} \alpha_i z_i^*  \right) = \Phi \left(\mathbf{a}\mathbf{z}^{*^T}\right)
\end{equation}

\noindent The MPP $\mathbf{z}^*$, reliability index $\beta = \left\lVert\mathbf{z}^*\right\rVert$ and the direction of gradient at MPP, $\mathbf{a}$ are therefore related to each other as

\begin{equation} \label{eq:31}
    \frac{\mathbf{z}^*}{\beta} = -\mathbf{a}
\end{equation}

\noindent Substituting \eref{eq:31} in \eref{eq:30}, results in  

\begin{equation} \label{eq:6228}
    P_f \approx P(L(\mathbf{z})<0)=\Phi (\mathbf{a}\mathbf{z}^{*^T}) = \Phi (-\beta \mathbf{a}\mathbf{a}^T) = \Phi(-\beta)
\end{equation}

\noindent Finally, the reliability of the system is given by 

\begin{equation}
    R = 1 - P_f = 1 - \Phi(-\beta) = \Phi(\beta)
\end{equation}

\subsection{Search for MPP} \label{MPP}
\begin{figure}[htbp!]  
\centering
\includegraphics[scale=0.4]{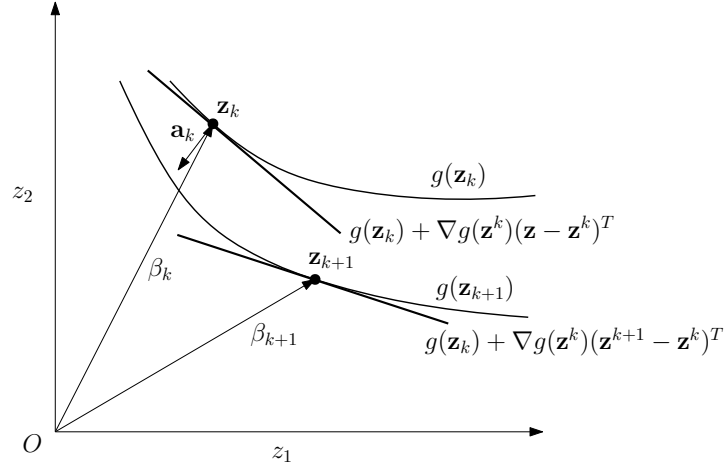}
\captionsetup{justification=centering}
\caption{Schematic representation of a MPP Search}
\label{fig:MPP_Search}
\end{figure}

Before proceeding with FORM/SORM, it has been made clear in Section \ref{FORM} that it is key to locate the MPP in the standard normal space. In this section, the commonly used MPP search algorithm will be discussed. The search is based on a recursive formula and on the linearization of the performance function, as highlighted in \fref{fig:MPP_Search} and explained below. 

Let MPP in the k'th iteration be $\mathbf{z}^k$. The performance function is expressed in a linearized form at $\mathbf{z}^k$ as shown by the lower line in \fref{fig:MPP_Search}. 

\begin{equation}
    g(\mathbf{z}) = g(\mathbf{z}^k)+\nabla g(\mathbf{z}^k)(\mathbf{z} - \mathbf{z}^k)^T
\end{equation}

\noindent Letting the linearized function to be zero, then the MPP $\mathbf{z}^{k+1}$ in the next iteration should be on the line as

\begin{equation} \label{eq:35}
    g(\mathbf{z}^{k+1}) = g(\mathbf{z}^k) + \nabla g (\mathbf{z}^k) (\mathbf{z}^{k+1} - \mathbf{z}^k)^T = 0
\end{equation}

\noindent From \eref{eq:31}, 

\begin{equation} \label{eq:36}
\begin{split}
    \mathbf{z}^k &= -\beta^k \mathbf{a}^k  \\ 
    \mathbf{z}^{k+1} &= -\beta^{k+1} \mathbf{a}^k
\end{split}
\end{equation}

\noindent Upon substituting \eref{eq:36} into \eref{eq:35}, one gets the followings: 

\begin{equation} \label{eq:37}
    \begin{split}
        g(\mathbf{z}^k)+\nabla g(\mathbf{z}^k)(\mathbf{a}^k)^T (\beta^k - \beta^{k+1}) &= g(\mathbf{z}^k) + \left\lVert \nabla g(\mathbf{z}^k) \right\rVert (\beta^k - \beta^{k+1}) = 0 \\
        \beta^{k+1} &= \beta^k + \frac{g(\mathbf{z}^k)}{ \left\lVert \nabla g(\mathbf{z}^k) \right\rVert} \\ 
        \mathbf{z}^{k+1} &= -\mathbf{a}^k \left( \beta^k + \frac{g(\mathbf{z}^k)}{ \left\lVert \nabla g(\mathbf{z}^k) \right\rVert} \right)
    \end{split}
\end{equation}

\noindent The last equation in \eref{eq:37} gives the recursive formula in estimating the MPP. As a starting point for the same, $\mathbf{z}^0 = 0$ is chosen. For the stopping criterion, the absolute change in the position of the iteration points is checked against a user defined criteria: $\left\lVert \mathbf{z}^{k+1} - \mathbf{z}^k \right\rVert \leq \epsilon$. 

\subsection{Second Order Reliability Method}\label{SORM}
The conventional Second Order Reliability Method (SORM) is based on second-order Taylor series approximation of the nonlinear limit state function so as to better approximate its curvature. The expansion about the MPP $\mathbf{z}^*$ is: 

\begin{equation}  \label{eq:38}
    g(\mathbf{z}) \approx g(\mathbf{z}^*) + \nabla g (\mathbf{z}^*) (\mathbf{z} - \mathbf{z}^*)^T + \frac{1}{2} (\mathbf{z} - \mathbf{z}^*) \mathbf{H}(\mathbf{z}^*) (\mathbf{z} - \mathbf{z}^*)^T
\end{equation}

\noindent where $\mathbf{H}(\mathbf{z}^*)$ is the Hessian matrix evaluated at MPP: 

\begin{equation}
\mathbf{H}(\mathbf{z}^*) = 
\begin{bmatrix}
{\frac{\partial^2 g}{\partial z_1^2}}\big\rvert_{\mathbf{z}^*} &  {\frac{\partial^2 g}{\partial z_1 z_2}}\big\rvert_{\mathbf{z}^*} & \cdots & {\frac{\partial^2 g}{\partial z_1 z_n}}\big\rvert_{\mathbf{z}^*} \\
{\frac{\partial^2 g}{\partial z_2 z_1}}\big\rvert_{\mathbf{z}^*} & {\frac{\partial^2 g}{\partial z_2^2}}\big\rvert_{\mathbf{z}^*} & \cdots & {\frac{\partial^2 g}{\partial z_2 z_n}}\big\rvert_{\mathbf{z}^*} \\
\vdots & \vdots & \ddots & \vdots\\
{\frac{\partial^2 g}{\partial z_n z_1}}\big\rvert_{\mathbf{z}^*} & {\frac{\partial^2 g}{\partial z_n z_2}}\big\rvert_{\mathbf{z}^*} & \cdots & {\frac{\partial^2 g}{\partial z_n^2}}\big\rvert_{\mathbf{z}^*} \end{bmatrix}
\label{eqn:39}
\end{equation}

\begin{figure}[htbp!]  
\centering
\includegraphics[scale=0.4]{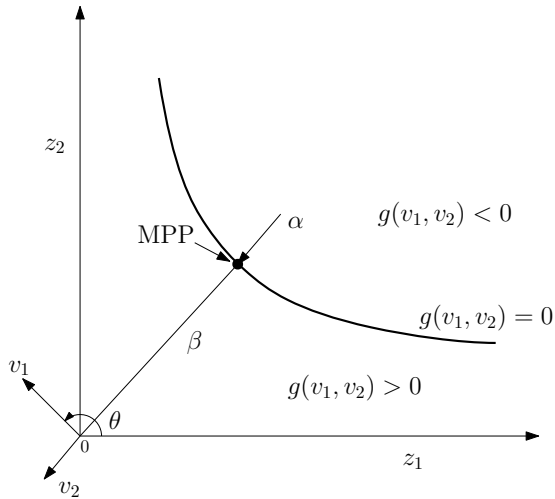}
\captionsetup{justification=centering}
\caption{Rotation of coordinates for SORM}
\label{fig:SORM_rot}
\end{figure}

\noindent A simple closed-form solution for the failure probability computation is derived using the theory of asymptotic approximation in \cite{breitung1984asymptotic} as: 

\begin{equation}
    P_f = \Phi(-\beta) \prod_{i=1}^{n-1} (1+\beta \kappa_i)^{-0.5}
\end{equation}

\noindent where $\beta$ is the reliability index obtained using FORM and $\kappa_i$ are the principal curvatures of the limit state at MPP. The procedure to obtain the first $(n-1)$ curvatures is detailed next. First, the standard Gaussian random variables $\mathbf{z}$ are rotated to another set of coordinates, $\mathbf{v}$, in such a way that the last component of the new set $v_n$, coincides with the $\alpha$ measure, which is the unit gradient vector of the limit state at the design point (see \fref{fig:SORM_rot} for the case of 2 random variables). This orthogonal transformation can be written as: 

\begin{equation}
    \mathbf{v}= \mathbf{R}\mathbf{z}
\end{equation}

\noindent where $\mathbf{R}$ is the rotation matrix. For the case of just 2 random variables, the rotation matrix will be 

\begin{equation}
    \mathbf{R} = 
    \begin{bmatrix}
    cos(\theta) & sin(\theta) \\
    -sin(\theta) & cos(\theta)
    \end{bmatrix}
\end{equation}

\noindent where $\theta$ is the rotation angle about the origin of the random space, as highlighted in \fref{fig:SORM_rot}. 

\noindent For the case of more than 2 random variables, the transformation matrix $\mathbf{R}$ is computed in two step:

\begin{steps}
  \item Construct the $\mathbf{R}_0$ matrix as follows: 
  \begin{equation}
      \mathbf{R}_0 = 
      \begin{bmatrix}
            1 & 0 & \cdots & 0 \\
            0 & 1 & \cdots & 0 \\ 
            \vdots & \vdots & \ddots & \vdots\\
            \alpha_1 & \alpha_2 & \cdots & \alpha_n 
      \end{bmatrix}
  \end{equation}
  \noindent where $\alpha_1,\alpha_2,\cdots,\alpha_n$ are the components of the unit gradient vector $\mathbf{a}$ at the design point, MPP. 
    \item Once the matrix $\mathbf{R}_0$ is computed, the Gram-Schmidt orthogonalization procedure is applied to $\mathbf{R}_0$ whose rows are $\mathbf{r}_{0_1},\mathbf{r}_{0_2},\cdots,\mathbf{r}_{0_n}$, to get the matrix $\mathbf{R}$ whose rows are $\mathbf{r}_1,\mathbf{r}_2,\cdots,\mathbf{r}_n$. The procedure for the above may be written as follows. The $n$'th row of $\mathbf{R}$ is simply $\mathbf{r}_n = \mathbf{r}_{0_n}$, while the other rows are computed in a backward order using the formula given below: 
    \begin{equation}
        \mathbf{r}_{k} = \mathbf{r}_{0_k} - \sum_{j=k+1}^{n} \left[ \frac{\mathbf{r}_j \mathbf{r}_{0_k}^T}{\mathbf{r}_j \mathbf{r}_j^T} \mathbf{r}_j    \right]
    \end{equation}
\end{steps}

Once the $\mathbf{R}$ matrix is obtained, a new matrix $\mathbf{A}$, whose elements are denoted by $A_{ij}$ is computed as follows:

\begin{equation}
A_{ij} = \frac{\left(\mathbf{R} \mathbf{H} \mathbf{R}^T     \right)_{ij}}{|\nabla g(\mathbf{v}^*)|}
\end{equation}

\noindent where $H$ is the Hessian matrix evaluated at the design point. Since $v_n$ coincides with the $\beta$-vector, the last row and column of the $\mathbf{A}$ matrix and the last row in the $\mathbf{v}$ vector are dropped. Finally the limit state function can be written using a second-order approximation in the rotated space $\mathbf{v}$ as: 

\begin{equation}
    \mathbf{Y}_n = \beta + \frac{1}{2} \mathbf{Y}^T A \mathbf{Y}
\end{equation}

\noindent where $\mathbf{A}$ is a matrix of size $(n-1)\times(n-1)$. Moreover, the curvatures $\kappa_i$ required for computing the probability of failure are just the eigenvalues of this reduced $\mathbf{A}$ matrix.

\subsubsection{Limitation of conventional FORM/SORM for reliability analysis of VSCL composites with cutouts}

As it can be observed from \eref{eq:dGdzi} and all throughout the derivation for FORM/ SORM in the earlier section, one needs to have an explicit representation of all input random variables $z_i$ within the definition of the limit state function $g(\mathbf{z})$. In our chosen particular problem for reliability analysis of composite plates with cutouts and including variability in geometric parameters like the ply thickness, it becomes unmistakably clear that the conventional approach with FORM/SORM is not the one to be advocated for this particular application. For example, by looking at \eref{eq:ABDmatrix}, one can clearly infer that, in case of random ply thickness, the integral limits become random, which makes the process of deriving an explicit representation of limit state function, $g(\mathbf{z})$ to be nearly intractable. All the more, by looking at \fref{fig:XFEM_mesh}, it is would be impractical to derive the gradient of the limit state function w.r.t the limit state function whose zero values defines the cutout contour. Hence, there exists a clear need to come up with a surrogate model that approximates the limit state function, whose discussion will be taken up in the following section.

\section{ANN based surrogate model for SRA} \label{sec:4}
In the earlier section, we had highlighted the need for surrogate modelling of limit state function for the problem in hand. Some of the most commonly used surrogate models for SRA are the Response Surface Method (RSM), Kriging \cite{romero2004construction,kaymaz2005application,bichon2008efficient}, Moment method, Neural Network based SRA, to name a few.

One of the first and primal works in applying Artificial Neural Networks (ANN) in SRA was by Hornik\textit{.et.al} \cite{HORNIK1989359}. For structural reliability analysis purposes, the two most commonly used ANN architectures are - back propagation multi-layer network and radial basis functions (RBF) \cite{CHOJACZYK201578}. The reader may refer to the works of Hurtado\textit{.et.al}\cite{hurtado2001neural} for an extensive state of the art discussion on general application of neural networks in stochastic mechanics. Once an adequate surrogate model is defined between the random input variables and its corresponding stochastic output using ANN, it is necessary to define a methodology wherein the probability of failure, $P_f$, can be computed using this ANN approximated meta-model. One methodology is to run an ANN-based Monte Carlo simulations for computing $P_f$, wherein only a finite number of finite element calculations are required for building the data set needed for ANN training process. Gradient-based methods like ANN coupled with FORM/SORM are found to abate the convergence issues normally encountered in conventional FORM/SORM (c.f. Section \ref{ConvFORMSORM}), especially when dealing with highly non-linear limit state function. 

In section \ref{Sec5}, the reader is introduced to the basics of a neural network architecture, followed by the procedure to derive the first and second order derivative of limit state function w.r.t each of the input random parameters.

\subsection{Overview of an Artificial Neural Network}\label{Sec5}

\begin{figure}[htbp!]  
\centering
\includegraphics[width=1\textwidth]{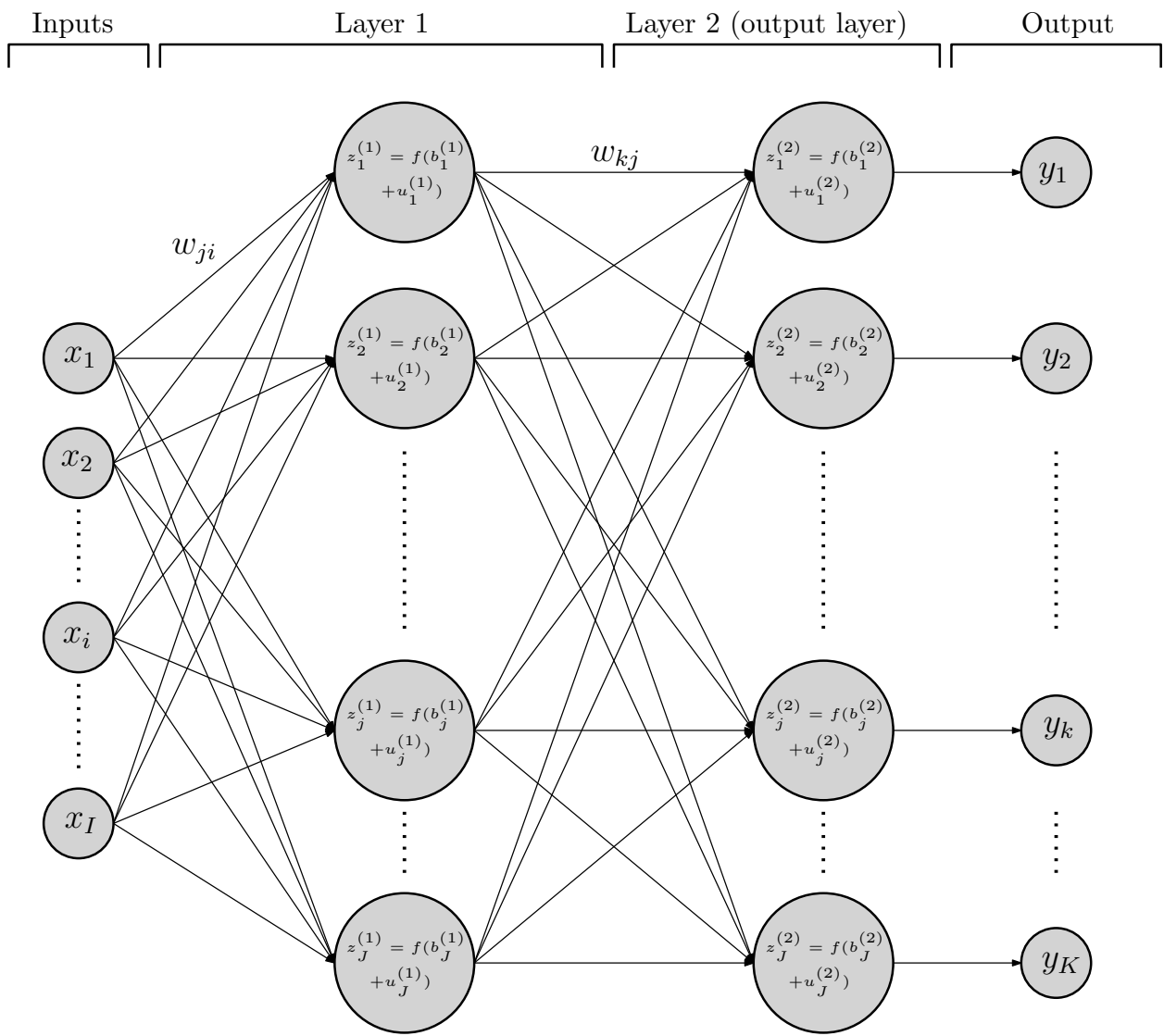}
\captionsetup{justification=centering}
\caption{Feedforward neural network with one hidden layer case.}
\label{fig:ANN_sketch}
\end{figure}

In an Artificial Neural Network (ANN), each layer containing neurons gets connected to the next layer. \fref{fig:ANN_sketch} shows a schematic representation of a feedforward neural network with a two-layer perceptron. As it can be seen from the outline sketch, a neuron in an ANN is considered to be a processing element. For the case of Layer 1, the output $z_j$ from the $j^{\rm th}$ neuron is calculated as follows: first, all its inputs from the neurons in the preceding layer $x_i$ are multiplied with their respective connection weights $w_{ji}$ and added up as shown below. 

\begin{equation} \label{eq:3_1}
    u_j = \sum_{i=1}^{I} w_{ji}x_i
\end{equation}

Before passing this weighted sum to the activation function, a small bias $b_i$ is added to the same. Finally, an activation function $f(\bullet)$ is needed to transform this biased weighted sum of the inputs $u_j$ into the output signal of the $j^{\rm th}$ neuron in Layer 1 as shown below. 

\begin{equation} \label{eq:221}
z_j^{1} = f \left( \sum_{i=1}^{I} w_{ji} x_i + b_j^{1} \right)
\end{equation}

\noindent where $I$ is the number of neurons in the input layer. Similar computations are carried out at each neuron within each layer of a feedforward neural-network until the output gets generated at the output layer.   

The optimal design of an ANN depends on a number of parameters. The setting up stage of an ANN can generally be divided into two stages: training and testing phase. ANN training phase is when the weight and biases of a user-defined ANN gets optimized/trained using a set of input-output numerically generated samples, with sufficient variation for each ANN input parameter from their respective distribution. In the testing phase, the accuracy of trained ANN in predicting the response to new unobserved data is checked for, so as to quantify the performance of a trained ANN. 
\noindent Several hyperparameters influence the performance of an ANN trained using back propagation. A detailed monologue for the same can be found in the works of Goodfellow.et.al \cite{Goodfellow-et-al-2016}. For the interest of the reader, a brief introduction to the same is provided below.   

\underline{Train/test ratio}: Typically, the ratio of training to test dataset size normally ranges between $80/20$ to $95/5$ ~\cite{Gulikers}. In this work, a training ratio of $80/20$ is chosen (meaning $80\%$ of samples are used for training the ANN, while the rest $20\%$ are used for testing).

\underline{Data normalization}: 
It is a good practice to employ data  normalization so as to avoid any unwanted biases towards select dimensions in the input-output space. Among several available pre-processing techniques like Min-Max, Z-score and Decimal Scaling Normalization to name a few, Min-Max normalization is used in current studies. 

\underline{ANN architecture}:
While the number of neurons in the input and output layer of an ANN is known \textit{a-priori}, the choice of number of hidden layers and the nodes per hidden layer remains open. It has been proven by \citet{Cybenko1989} that a one-layer ANN can sufficiently approximate any nonlinear function.  While we restrict ourselves to a shallow one hidden layer ANN architecture, the choice of number of neurons in the hidden layer is discussed in the numerical section \ref{Numerical_examples}.

\underline{Activation function}:
Choosing the right neuron activation function $f(\bullet)$ in \eref{eq:221} is an appropriate consideration, because it directly affects how the input data gets formatted. There exists several activation functions in the literature, depending on the range of normalized input data. In this study, given that we use a Min-Max data normalization wherein the input data gets normalized to within a range of [-1,1], hyperbolic tangent function is used as activation function for hidden layer neurons, while linear activation is used in the output layer neuron.  

\underline{Cost function}
Training a feedforward ANN is more or less an optimization problem, for which a cost function has to be defined. There exists several cost functions to choose from, of which the Mean Squared activation function is chosen in the present study. The mean squared cost function is widely used in statistics, mainly for regression problems. The weights and bias of an ANN are optimized using stochastic gradient descent algorithm ~\cite{Papalambros:2003:POD:1202943}. \newline
With the help of a trained ANN, the expressions to compute both the First Order Derivative (FOD)and the Second Order Derivative (SOD) of the ANN derived function approximation w.r.t each of the input parameters are presented next. The detailed derivation for the same can be found elsewhere ~\cite{deng2005structural,tawfik2018neural}. However, a limitation of the back-propagation algorithm is that different ANNs with the same architecture and training data set, gives rise to different difference. This is due to the gradient decent algorithm used to minimize the error function of an ANN, which gets trapped trapped on a different error local minima in each ANN \cite{oparaji2017robust}.

\subsection{ANN derived First Order Derivative}\label{ANN_FOD}

The first order derivatives can be computed by using the chain rule of differentiation. With reference to the ANN architecture as in \fref{fig:ANN_sketch}, the derivative of the ANN output $y_k$ w.r.t. to its input $x_i$ can be obtained as follows, using the assumption in our study that the transfer function in the output layer is linear and that in the hidden layer is hyperbolic tangential. Referring to \eref{eq:221} and \fref{fig:ANN_sketch}, the derivative of the output $y_{k}$ w.r.t the input $x_{i}$ can be obtained using the chain rule of differentiation as follows: 

\begin{align} \label{eq:331}
\frac{\partial y_k}{\partial x_i} &= \frac{\partial y_k}{\partial z_j} \frac{\partial z_j}{\partial x_i}, \qquad j=1,\cdots,J 
\end{align}

\noindent
where $J$ is the number of neurons in the hidden layer. 

\begin{align}
    \frac{\partial y_k}{\partial z_j} &= w_{kj}  \label{eq:32} \\
    \frac{\partial z_j}{\partial x_i} &= w_{ji} (1-z_j^2) \label{eq:33} 
\end{align}

\noindent The term that multiplies the weight factor in \eref{eq:33} is the derivative of the activation function (hyperbolic tangential). Substituting \eref{eq:32} and \eref{eq:33} into \eref{eq:331}, one gets the following relation for FOD: 

\begin{equation}
    \frac{\partial y_k}{\partial x_i} = \sum_{j=1}^{J} w_{kj} w_{ji} (1-z_j^2)
\end{equation}

\subsection{ANN derived Second Order Derivative }\label{ANN_SOD}

In a similar manner, the second order derivative of the ANN output w.r.t the random inputs can be estimated using the chain rule \cite{deng2005structural}: 

\begin{equation}
\begin{split}
\frac{\partial^2 y_k}{\partial x_j \partial x_{j'}} &= \frac{\partial}{\partial x_{j'}} w_{kj} w_{ji} (1 - z_j^2) \\ &= w_{kj} w_{ji} (-2z_j) \frac{\partial z_i}{\partial x_{j'}} \\
&= w_{kj} w_{ji} w_{j'i} (-2z_j) (1-z_j^2)  \qquad j,j'=1,\cdots,I
\end{split}
\end{equation}

\section{ANN based Probability of Failure estimation for Variable Stiffness Composite Laminate Plates} \label{ANN_VSCL_prob}

\begin{figure}[htbp!]  
\centering
\includegraphics[width=0.7\textwidth]{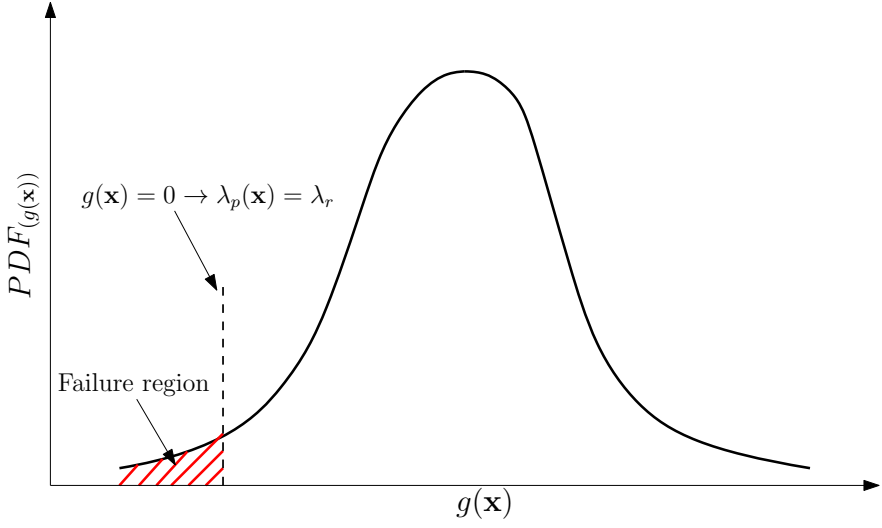}
\captionsetup{justification=centering}
\caption{Sample PDF plot of limit state function along with representation of failure region (g($\mathbf{x}$)$<$0)}
\label{fig:gX}
\end{figure}

With the aid of a brief overview over feed-forward ANN architecture and its functional derivatives, the emphasis in this section is placed on defining a suitable limit state function followed by its gradient interpretation, which will later be used in ANN based FORM-SORM. In this work, a limit state function, $g(\mathbf{x})$, is defined based on the fundamental free-vibration frequency of a VSCL plate, as shown below:   

\begin{equation} \label{eq:12}
g(\mathbf{x})=\frac{\lambda_p(\mathbf{x})}{\lambda_r}-1
\end{equation}

\noindent where $\lambda_p(\mathbf{x})$  corresponds to the random fundamental frequency computed for an input experimental point $\mathbf{x}$ of dimension $n$ and $\lambda_r$ stands for a user-defined critical value below which the structure is considered to be failed. \fref{fig:gX} shows a sample PDF plot of the above limit state function definition, where the probability of failure region is also highlighted. 

In ANN-based FORM-SORM, partial derivatives of the performance function $g(\mathbf{x})$ w.r.t each of the input random variables $x_i, i=1,\cdots,n$ are required and are obtained as follows: 

\begin{equation}
    \begin{split}
        \frac{\partial g(\mathbf{x})}{\partial x_i} &= \frac{\partial g}{\partial \lambda_p} \frac{\partial \lambda_p(\mathbf{x})}{\partial x_i} = \frac{1}{\lambda_r} \frac{\partial \lambda_p(\mathbf{x})}{\partial x_i} \approx \frac{1}{\lambda_r} \left(\frac{\partial \lambda_p(\mathbf{x})}{\partial x_i} \right)_{ANN} \qquad i=1,2,\cdots,n \\
        \frac{\partial^2 g(\mathbf{x})}{\partial x_i \partial x_j} &= \frac{1}{\lambda_r} \frac{\partial^2 \lambda_p(\mathbf{x})}{\partial x_i \partial x_j} \approx \frac{1}{\lambda_r} \left( \frac{\partial^2 \lambda_p(\mathbf{x})}{\partial x_i \partial x_j} \right)_{ANN}  \qquad i,j=1,2,\cdots,n
    \end{split}
\end{equation}
\noindent One can see from the above equation that the partial derivatives of the fundamental frequency of the system w.r.t each of the input random variables $x_i$ are required. The identities within the brackets $(\cdot)_{ANN}$ are the ones that we try to approximate using the trained ANN, as detailed in Section \ref{ANN_FOD} and \ref{ANN_SOD}. Moreover, one should keep in mind that as per the definition of FORM-SORM, the derivatives of the limit state function w.r.t input parameters in their standard normal distribution space $\mathcal{Z}$ is required (c.f. \eref{eq:dGdzi}). However, the limit state function derivatives that are obtained using the trained ANN is based on the original input  distribution of random parameters, say $\mathbf{x} \approx \mathcal{N}(\mu_i,\sigma_i^2);i=1,\cdots,n$. Hence, a transformation as defined in \eref{eqn:transformation} is required, coupled with \eref{eq:16}. 

\begin{equation} \label{eqn:transformation}
    \begin{split}
         \frac{\partial g}{\partial z_i} &= \frac{\partial g}{\partial x_i} \frac{\partial x_i}{\partial z_i} \approx \left( \frac{\partial g}{\partial x_i}  \right)_{ANN} \sigma_i   \\ 
         \frac{\partial^2 g}{\partial z_i \partial z_j}  &\approx \left( \frac{\partial^2 g}{\partial x_i  \partial x_j}\right)_{ANN}\sigma_i\sigma_j
    \end{split}
\end{equation}

\section{Sampling method for ANN training set generation} \label{ANN_sampling}

In ANN based reliability method, the accuracy of reliability estimates is strongly dependent on the quality of the samples used to train the network \cite{CHOJACZYK201578}. The general approach adopted for generating training data set is to randomly select the samples across the whole input parameter domain and the weights and the biases of the ANN gets iteratively adjusted until a required level of accuracy is achieved in level set function approximation ~\cite{PAPADRAKAKIS1996145,anjum1997response,pu2006application}. However, such traditional approach has been found to be less efficient since only a low fraction of the randomly generated samples lay in a region located close to $g(\mathbf{x})=0$ and that too next to the MPP. 

\subsection{Monte Carlo Simulation (MCS)} \label{Sec3}

The integral in \eref{eq:Pf} can be evaluated  by sampling large number of samples from the input parametric space as shown below:

\begin{equation} \label{eq:777}
    P_{f_{MCS}} = \frac{1}{N} \sum_{i=1}^{N} \mathbb{I}(g(\mathbf{x}_i)<0)
\end{equation}

\noindent where $\mathbb{I}(\bullet)$ stands for the indicator function and $N$ stands for the number of samples. The samples are chosen from the joint PDF $f_{X}(\mathbf{x})$. However, one of the major drawbacks of naive MCS estimator is that they tend to have an extremely high variance in $P_f$ prediction, especially when failure probability is large only on a small fraction of samples \cite{lopes2010reliability}. To study the convergence of failure probability estimates as a function of number of samples, a coefficient of variation measure for \eref{eq:777} is used as given below \cite{jian2017two}:

\begin{equation} \label{eq:1900}
    COV_{P_{f_{MCS}}} = \frac{\sqrt{\text{var}(P_{f_{MCS}}})}{P_{f_{MCS}}} = \sqrt{\frac{1-P_{f_{MCS}}}{N P_{f_{MCS}}}} \approx \frac{1}{\sqrt{N P_{f_{MCS}} }}
\end{equation}

For example, for a $10\%$ of coefficient of variation, one would require a minimum number of sample size of $N= 10^{k+2}$ for estimating a failure probability of order $10^{-k}$ \cite{dubourg2014meta}. 

\subsection{Importance Sampling and its adaptive variants}\label{IS}

To overcome this limitation, numerous variants have been developed, like the Subset Sampling (SS) and Importance Sampling (IS), to name a few. Importance sampling simulates from a biasing distribution that may be different from the true underlying joint distribution of input random parameters; however a correction is made for this mismatch by weighing the samples with an appropriate ratio of densities \cite{LiMarzouk}. By focusing samples in regions where the integrand in \eref{eq:Pf} is large, IS can reduce the variance of a Monte Carlo estimator of an integral \cite{Liu:2008:MCS:1571802}. In short, \eref{eq:Pf} gets modified as follows: 

\begin{equation} \label{eq:IS_1}
    P_{f_{IS}} = \int_{g(\mathbf{x})\leq 0} \frac{f_{X}(\mathbf{x})}{h_X(\mathbf{x})} h_X(\mathbf{x}) \rm{d}x
\end{equation} 

\noindent where $h_X(\mathbf{x})$ is the bias distribution or also termed as the instrumental density. The unbiased importance sampling based $P_f$ estimator then becomes  

\begin{equation}\label{equ:IS_2}
    P_{f_{IS}} = \frac{1}{N} \sum_{i=1}^{N} \mathbb{I}(g(x_i)<0) \frac{f_X(x_i)}{h_X(x_i)} 
\end{equation}

\noindent where $\{x_1,\cdots,x_N \}$ are set of samples drawn from the instrumental density $h_X(\mathbf{x})$. As per the central limit theory, the IS estimator is unbiased and and its quality can be predicted using the variance of estimation which is given below \cite{dubourg2013metamodel}: 

\begin{equation} \label{eq:1555}
    Var[P_{f_{IS}}] = \frac{1}{N-1}\left( \frac{1}{N} \sum_{i=1}^{N} \mathbb{I}(g(x_i)<0) \frac{f(x_i)^2}{h(x_i)^2} - P_{f_{IS}}^2 \right)
\end{equation}

The effectiveness of IS depends on the selection of an appropriate $h_X(\mathbf{x})$, also known as the instrumental density in literature, such that the probability sampling in \eref{equ:IS_2} can be prioritized onto the region of greatest importance. Au\textit{.et.al} \cite{au1999new} proposed the concept of sampling from the instrumental density as states of a Markov chain, which asymptotically converges to the optimal instrumental sampling desnity to be represented with increase in samples. Afterwards, a kernel sampling density was constructed using these samples which was then later used as $h_X(\mathbf{x})$ in IS. Dubourg\textit{.et.al} \cite{dubourg2014meta} proposed a meta-model based importance sampling method where the meta-model to approximate the limit state function is based on the kriging approach. However, the probabilistic nature of the meta-model was used to construct a quasi-optimal $h_X(\mathbf{x})$. 

\section{Global Reliability Sensitivity Analysis using Artificial Neural Network} \label{ANN_reliability}

Apart from predicting the failure probability of a structure, it is equally important as to how to reduce the failure probability by carrying out a total sensitivity studies. Several work have been published along these lines. Broadly there exists two ways to compute the reliability based sensitivity studies, namely the Parametric Reliability Sensitivity (PRS) and Global Reliability Sensitivity (GRS). For example, Dubourg.\textit{et.al.} \cite{dubourg2013metamodel} used the score function approach to derive the sensitivity of the failure probability w.r.t the parameter set. Assuming each of the input random variables are Gaussian, the parameter set for the problem consists of $\bm{\theta}_i = (\mu_i,\sigma_i); i = 1,\cdots,n$, where $\mu_i$ and $\sigma_i$ corresponds to the mean and standard deviation of the i'th input random variable. Recently, Zhang\textit{.et.al} \cite{zhang2019effective} proposed the concept of using the principle of maximum entropy (MaxEnt) to approximate the reliability based sensitivity index using a performance function approximated using the multiplicative dimensional reduction method. Pengfei \textit{.et.al} proposed a method to estimate both PRS and GRS, based on importance samples obtained using quasi-optimal importance sampling passed through an Adaptive Kriging approximated performance function. However, what we are interested here is to compute the total sensitivity indices and for that, we will be using the procedure outlined in \cite{wei2012efficient}. The reader is advised to refer to the work of Wei\textit{.et.al} \cite{wei2012efficient} for the detailed derivations and procedures to be followed for both MCS based and Importance Sampling based Total Sensitivity Index estimation for each input random parameters on the structural failure probability.




\section{Numerical Examples} \label{Numerical_examples}

In this section, the proposed ANN based methodology for predicting the failure probability of a VSCL plate in the presence of both geometric and material uncertainties is studied. An in-house code was developed to numerically estimate the first natural mode of vibration of a VSCL plate with an explicit cutout representation using the XFEM, given a deterministic set of input parameters. The influence of various geometric parameters, viz., plate aspect ratio, $a/b$, ply thickness $t_k$, fiber orientation angle $<\Theta_0^k | \Theta_1^k>$ and the geometry of the cutout on the probability of the failure is systematically studied. The plate is assumed to be simply supported on all edges and is discretized using four noded shear flexible structured quadrilateral element with five dofs per node $(u_o,v_o,w_o,\beta_x,\beta_y)$. All computations were carried out on a desktop with Intel Quad Core i7, 2.40 GHz CPU with 8 GB RAM. 

\subsection{Validation of deterministic code with literature results}

Before proceeding with a systematic parametric study, the deterministic results from the present work is validated against the results in the literature for a square VSCL plate (a=b=1m) without cutouts. All parameters are assumed to be deterministic and the plate is made up of three layers with the following layup: $\left[<30 \degree,0\degree>,<45\degree,90\degree>, <30\degree,0\degree>\right]$. Total thickness of the plate is taken as 0.01 m. The material properties are: $E_1=$ 173 Gpa, $E_2=$ 7.2 GPa, $G_{12}=G_{13}= G_{23}$ = 3.76 GPa, $\nu_{12}=\nu_{13}=\nu_{23}=$ 0.29 and $\rho=$ 1540 kg/m$^3$. Table \ref{tab:7} presents the convergence of the first five fundamental frequencies with mesh refinement. Based on a progressive refinement, a structured mesh of 30$\times$30 is found to be adequate to model the plate.

\begin{table}[htbp!]
\captionsetup{justification=centering}
\caption{Natural frequency, $\omega$ $(rad/s)$ for a square simply supported VSCL plate without a cutout}
\label{tab:7}       
\centering
\begin{tabular}{ccccc}
\hline\noalign{\smallskip}
Ref.~\cite{ribeiro2014review} & \multicolumn{4}{c}{Mesh size (Present work)} \\ \cmidrule{2-5}
 & $10\times 10$  & $20 \times 20$ & $30 \times 30$ &  $\% \Delta_{30 \times 30}$ \\
\noalign{\smallskip}\hline\noalign{\smallskip}
309.1 & 315.0 & 309.9 & 309.1  & 0.0136  \\
503.3 & 530.5 & 509.3 & 505.9  & 0.5245  \\
852.1 & 954.0 & 867.4 & 854.3  & 0.2627  \\
1143.5 & 1203.4 & 1144.0  & 1134.3 & 0.8053 \\
1297.3 & 1416.1 & 1324.8 & 1296.2 & 0.0805  \\
\noalign{\smallskip}\hline
\end{tabular}
\end{table}

Next, the probability of failure predicted using the ANN trained FORM/SORM method is validated with the results available in the literature for Constant Stiffness Composite Laminate (CSCL)~\cite{tawfik2018neural} in the absence of a cutout. The plate is assumed to be a square with thickness $h=$ 0.1. Table \ref{tab:Stat_distr_mech_prop_valid} summarizes the statistical distribution of all the input variables. All units are in SI. A stacking sequence of $[0 \degree /45\degree/-45\degree/90\degree]$ is considered. A 1000 samples were generated using the statistical distribution of the input parameters given in Table \ref{tab:Stat_distr_mech_prop_valid}. Out of these 1000 samples, 500 were used for training the ANN with an architecture of [11-10-1] (similar to \cite{tawfik2018neural}), while the remaining 500 samples were shared between testing and validation phase. Hyperbolic tangent activation function was used for the hidden layer (Logistic Sigmoid activation function was used for the hidden layer in Tawfik.et.al \cite{tawfik2018neural}), except for the output layer where a linear activation function is used and $\omega_r/\omega_p = 0.97$ is used in this study.
 
\begin{table}[htbp!]
\captionsetup{justification=centering}
    \caption{Statistical distribution of basic random variables used for validation of NN based FORM/SORM Matlab code validation \cite{tawfik2018neural}.}
    \centering
    \begin{tabular}{l c c c c c c c c}
    \hline\noalign{\smallskip}
         Property & $E_{11}$ & $E_{22}$ & $G_{12}$ & $G_{23}$ & $G_{13}$ & $\nu_{12}$ & $\rho$ & $\Delta \theta$ \\
         \hline\noalign{\smallskip}
         Mean & 16.48 & 1.4 & 0.87 & 0.45 & 0.87 & 0.334 & 1000 & $0 \degree$ \\
         Std Dev & 0.61 & 0.05 & 0.052 & 0.014 & 0.052 & 0.01 & 36 & 1.8 $\degree$ \\ 
         \hline\noalign{\smallskip}
    \end{tabular}
    \label{tab:Stat_distr_mech_prop_valid}
\end{table}
 
\begin{figure}[htpb!] 
\newlength\figureheight
\newlength\figurewidth
\setlength\figureheight{8cm} 
\setlength\figurewidth{8cm}
\centering 
%
%
\begin{tikzpicture}

\begin{axis}[%
width=0.95092\figurewidth,
height=\figureheight,
at={(0\figurewidth,0\figureheight)},
scale only axis,
xmin=0,
xmax=120,
xlabel={Iteration number},
ymode=log,
ymin=1e-09,
ymax=0.01,
yminorticks=true,
ylabel={Mean Squared Error},
legend style={legend cell align=left,align=left,draw=white!15!black}
]
\addplot [color=red,solid,line width=2.0pt]
  table[row sep=crcr]{%
0	0.00416199985261351\\
1	0.000657928481571665\\
2	0.000110943552023313\\
3	1.59271073040822e-05\\
4	1.13891163735528e-05\\
5	7.64595075564836e-06\\
6	7.55478466995303e-07\\
7	2.07792855833952e-07\\
8	1.56406185543558e-07\\
9	1.40174783475285e-07\\
10	1.32162632811795e-07\\
11	1.27319258817953e-07\\
12	1.23971170082803e-07\\
13	1.2142553735451e-07\\
14	1.19354728077226e-07\\
15	1.17586493819239e-07\\
16	1.16021696832637e-07\\
17	1.14598803606705e-07\\
18	1.13277221014056e-07\\
19	1.12028863083249e-07\\
20	1.10833589389301e-07\\
21	1.09676581734275e-07\\
22	1.08546736361184e-07\\
23	1.07435602490031e-07\\
24	1.06336618042782e-07\\
25	1.05244509719781e-07\\
26	1.04154791712665e-07\\
27	1.03063338524637e-07\\
28	1.01966032071637e-07\\
29	1.01631985112995e-07\\
30	8.87012252996219e-08\\
31	7.47745248182039e-08\\
32	6.19483753575462e-08\\
33	5.10772464348525e-08\\
34	4.31562928518418e-08\\
35	3.81095943887221e-08\\
36	3.45795443148449e-08\\
37	3.18807091205667e-08\\
38	2.98506261177981e-08\\
39	2.83663279485658e-08\\
40	2.72821865116822e-08\\
41	2.64712753960607e-08\\
42	2.5843323923345e-08\\
43	2.53402556643727e-08\\
44	2.49258309448052e-08\\
45	2.4577029162948e-08\\
46	2.42785367114274e-08\\
47	2.40195821979149e-08\\
48	2.37921820217524e-08\\
49	2.35901619132738e-08\\
50	2.34085944496885e-08\\
51	2.32434599768861e-08\\
52	2.30914288694913e-08\\
53	2.29497107141388e-08\\
54	2.28159413616625e-08\\
55	2.26880924506155e-08\\
56	2.25643953011272e-08\\
57	2.24432748430257e-08\\
58	2.232329106326e-08\\
59	2.22030862535344e-08\\
60	2.20813366351426e-08\\
61	2.19567070220122e-08\\
62	2.18278072027039e-08\\
63	2.16931487494284e-08\\
64	2.15511010345223e-08\\
65	2.13998453825768e-08\\
66	2.12373265505515e-08\\
67	2.10612011746444e-08\\
68	2.10094594470199e-08\\
69	1.6701102948832e-08\\
70	1.24748321334302e-08\\
71	9.62069152262993e-09\\
72	8.24814037802685e-09\\
73	7.06955857625591e-09\\
74	6.37126250194454e-09\\
75	5.91518391478047e-09\\
76	5.57167198466436e-09\\
77	5.29141312706025e-09\\
78	5.04918805193324e-09\\
79	4.83028888256699e-09\\
80	4.62490952673528e-09\\
81	4.42589009966071e-09\\
82	4.2276808775069e-09\\
83	4.02615959467955e-09\\
84	3.81910819309301e-09\\
85	3.60717878494866e-09\\
86	3.39475573860893e-09\\
87	3.18959113323581e-09\\
88	3.00042696654762e-09\\
89	2.83363989702065e-09\\
90	2.69132024936824e-09\\
91	2.57184610247137e-09\\
92	2.47175737662598e-09\\
93	2.38734281223867e-09\\
94	2.31535899492923e-09\\
95	2.25315521295025e-09\\
96	2.19861256012146e-09\\
97	2.1500702330304e-09\\
98	2.10625484902013e-09\\
99	2.06620701071408e-09\\
100	2.02921058568603e-09\\
101	1.99473146247326e-09\\
102	1.96236805060932e-09\\
103	1.93181258831756e-09\\
104	1.9028215344287e-09\\
105	1.87519368247381e-09\\
};
\addlegendentry{Training error};
\addplot [color=blue,dashed,line width=2.5pt, mark=o]
  table[row sep=crcr]{%
0	0.00346264620452094\\
1	0.000674010634043657\\
2	0.000113633024602159\\
3	1.44576270756996e-05\\
4	1.11583751252553e-05\\
5	7.48202984482026e-06\\
6	7.98693352364806e-07\\
7	2.69449561616907e-07\\
8	2.16770296672356e-07\\
9	1.99819429129495e-07\\
10	1.9068066337267e-07\\
11	1.84326836554396e-07\\
12	1.79246673488319e-07\\
13	1.74884184478457e-07\\
14	1.70987329958084e-07\\
15	1.67413904753285e-07\\
16	1.64069580635064e-07\\
17	1.60886337645253e-07\\
18	1.57813685418954e-07\\
19	1.54814257434318e-07\\
20	1.51861134707243e-07\\
21	1.48935978143976e-07\\
22	1.46027570210215e-07\\
23	1.43130544708769e-07\\
24	1.40244170967303e-07\\
25	1.37371126979141e-07\\
26	1.34516260851617e-07\\
27	1.31685398083543e-07\\
28	1.28884293612203e-07\\
29	1.30634090022193e-07\\
30	1.09976917863334e-07\\
31	9.07217722873183e-08\\
32	7.65526910213278e-08\\
33	6.6916333052033e-08\\
34	6.06447980683881e-08\\
35	5.65897670838048e-08\\
36	5.35770413247127e-08\\
37	5.11883305895974e-08\\
38	4.93261843037017e-08\\
39	4.78939570131181e-08\\
40	4.67827154480978e-08\\
41	4.59019947477756e-08\\
42	4.51861459358609e-08\\
43	4.45898522068923e-08\\
44	4.40822578040451e-08\\
45	4.36421589519989e-08\\
46	4.32546597425204e-08\\
47	4.29090050532328e-08\\
48	4.25971976561598e-08\\
49	4.23131089575281e-08\\
50	4.20519003771494e-08\\
51	4.18096420854704e-08\\
52	4.15830572039051e-08\\
53	4.13693451612338e-08\\
54	4.11660545262833e-08\\
55	4.09709865484325e-08\\
56	4.07821176171035e-08\\
57	4.05975331242094e-08\\
58	4.04153676921138e-08\\
59	4.02337480771702e-08\\
60	4.00507357242922e-08\\
61	3.98642662184167e-08\\
62	3.96720829303619e-08\\
63	3.94716621012289e-08\\
64	3.92601265436646e-08\\
65	3.90341451655916e-08\\
66	3.87898157976545e-08\\
67	3.85225295752251e-08\\
68	3.78191414504793e-08\\
69	3.15293324617222e-08\\
70	2.32035357398362e-08\\
71	1.64002010713978e-08\\
72	1.31402083089339e-08\\
73	1.13695625434295e-08\\
74	1.04313258212363e-08\\
75	9.84558120694355e-09\\
76	9.41101194809018e-09\\
77	9.04996031463671e-09\\
78	8.72501665335672e-09\\
79	8.41592115291561e-09\\
80	8.11076066700521e-09\\
81	7.80237153763578e-09\\
82	7.48705532227377e-09\\
83	7.16458517907107e-09\\
84	6.83886701413513e-09\\
85	6.51836006952852e-09\\
86	6.2149264543036e-09\\
87	5.9402534901309e-09\\
88	5.70127000876215e-09\\
89	5.49798952547183e-09\\
90	5.32530109675097e-09\\
91	5.17640215679372e-09\\
92	5.04520124273332e-09\\
93	4.92716595355293e-09\\
94	4.81921310531688e-09\\
95	4.71926748579057e-09\\
96	4.62587072828962e-09\\
97	4.53794302222298e-09\\
98	4.45465386337694e-09\\
99	4.37534669405831e-09\\
100	4.29948885751785e-09\\
101	4.22663577059716e-09\\
102	4.15640529306866e-09\\
103	4.08846091669535e-09\\
104	4.0225030733436e-09\\
105	3.95826721668399e-09\\
};
\addlegendentry{Validation error};

\addplot [color=black,dotted,line width=1.5pt]
  table[row sep=crcr]{%
0	0.0040015039606651\\
1	0.000703500873611495\\
2	0.00012108560919006\\
3	1.92081314360596e-05\\
4	1.31751632775659e-05\\
5	8.67287395253467e-06\\
6	1.09077401613615e-06\\
7	3.93852648465386e-07\\
8	3.05975692497181e-07\\
9	2.74805531218438e-07\\
10	2.57747939731326e-07\\
11	2.45982448183105e-07\\
12	2.36797249263263e-07\\
13	2.29179089462007e-07\\
14	2.22652250492477e-07\\
15	2.1693683095169e-07\\
16	2.1184111968992e-07\\
17	2.07223496876456e-07\\
18	2.02975901010056e-07\\
19	1.99014788544387e-07\\
20	1.95275300058106e-07\\
21	1.91707143556535e-07\\
22	1.88271528032459e-07\\
23	1.84938762960772e-07\\
24	1.81686277161577e-07\\
25	1.78496907264458e-07\\
26	1.75357386015434e-07\\
27	1.72257025704479e-07\\
28	1.6918663742479e-07\\
29	1.68924598934739e-07\\
30	1.44119575401354e-07\\
31	1.20300943178465e-07\\
32	1.02432692835001e-07\\
33	9.05143780363165e-08\\
34	8.32571657530163e-08\\
35	7.92904823637459e-08\\
36	7.69073611228271e-08\\
37	7.51675711536024e-08\\
38	7.37289644933067e-08\\
39	7.24145192401836e-08\\
40	7.11476098604861e-08\\
41	6.99218280341326e-08\\
42	6.87554254191079e-08\\
43	6.76649803001566e-08\\
44	6.66576922387955e-08\\
45	6.57322945342683e-08\\
46	6.48822830021127e-08\\
47	6.40988046045832e-08\\
48	6.33725337237377e-08\\
49	6.26946542168062e-08\\
50	6.20572503200087e-08\\
51	6.14533688514074e-08\\
52	6.08769264762586e-08\\
53	6.03225606704142e-08\\
54	5.97854733547383e-08\\
55	5.92612876137726e-08\\
56	5.87459232443619e-08\\
57	5.82354902525725e-08\\
58	5.77261969963862e-08\\
59	5.72142692414132e-08\\
60	5.66958767925391e-08\\
61	5.61670650409446e-08\\
62	5.56236895084668e-08\\
63	5.50613522345305e-08\\
64	5.447533967749e-08\\
65	5.38605627806408e-08\\
66	5.32115011115182e-08\\
67	5.25221546880427e-08\\
68	4.80459046736833e-08\\
69	4.00717067376322e-08\\
70	3.22082517473974e-08\\
71	2.77076701760619e-08\\
72	2.51148575233973e-08\\
73	2.25909235217932e-08\\
74	2.09830807621427e-08\\
75	1.98368919003762e-08\\
76	1.89394722304999e-08\\
77	1.81934165626379e-08\\
78	1.75459142142886e-08\\
79	1.69635724974317e-08\\
80	1.64241597902413e-08\\
81	1.59118601672735e-08\\
82	1.54152034611431e-08\\
83	1.49264530206512e-08\\
84	1.44422714085e-08\\
85	1.39650270658084e-08\\
86	1.35034339727609e-08\\
87	1.30706302525497e-08\\
88	1.2679163041336e-08\\
89	1.23357515652019e-08\\
90	1.20398273666311e-08\\
91	1.17860409963657e-08\\
92	1.15675659551203e-08\\
93	1.13780301781424e-08\\
94	1.12120650417604e-08\\
95	1.10652266577346e-08\\
96	1.09338238003807e-08\\
97	1.08148035658847e-08\\
98	1.07056683264074e-08\\
99	1.06043921599814e-08\\
100	1.05093349397715e-08\\
101	1.04191639018412e-08\\
102	1.03327902927798e-08\\
103	1.0249323339074e-08\\
104	1.01680392174426e-08\\
105	1.00883596012021e-08\\
};
\addlegendentry{Testing error};

\end{axis}
\end{tikzpicture}%
\captionsetup{justification=centering}
\caption{Training, Validation and Testing errors of ANN. Case for the reliability estimates w.r.t the first mode of vibration studies for a SSSS $[0\degree/45\degree/-45\degree/90\degree]$ composite plate using a [11-10-1] ANN architecture with randomness in input parameters as given in table \ref{tab:Stat_distr_mech_prop_valid}. }
\label{fig:TVT_error}
\end{figure}
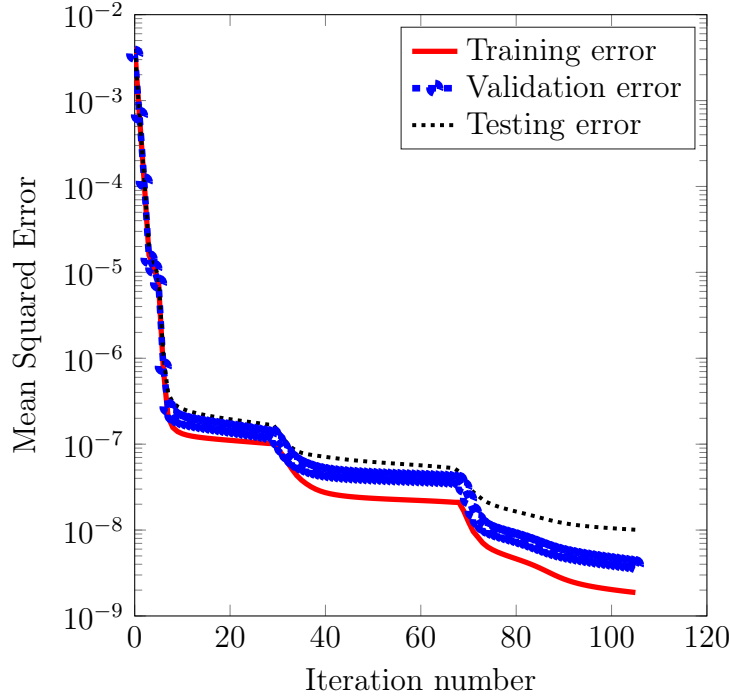

Tables \ref{tab:MPP_validation_Tawfik} - \ref{tab:Pf_beta_comparison} compares the predicted \textit{Most Probable Point} (MPP) for failure, the probability of failure $P_f$ and the reliability index $\beta$ using our in-house built MATLAB code with the results available in \cite{tawfik2018neural}. It is evident that a good agreement in predictions is achieved. The small differences in the magnitude can be attributed to: (a) the difference in the activation function selection for the hidden layer and (b) the choice of the plate deformation theory (FSDT is used in the current work, whilst, the reference work \cite{tawfik2018neural} is based on Third order Shear Deformation Theory (TSDT)). The training, validation and testing errors are plotted against the iteration number in \fref{fig:TVT_error}. After nearly 105 iterations, the training, validation and testing error were found to be 1.86$\times$10$^{-9}$, 3.96$\times$10$^{-9}$ and 1.01$\times$10$^{-9}$ respectively. 

\begin{table}[htbp!]
\captionsetup{justification=centering}
    \caption{MPP for a square simply supported composite laminate with $[0\degree/45\degree/-45\degree/90\degree]$ layup.}
    \centering
    \begin{tabular}{lrr}
    \hline\noalign{\smallskip}
         Random variables & Present Work & Tawfik.et.al \cite{tawfik2018neural} \\
         \hline\noalign{\smallskip}
         $E_{11}$ (GPa)  & 16.098 & 16.099 \\
         $E_{22}$ (GPa)  & 1.3903 & 1.3828 \\
         $G_{12}$ (GPa)  & 0.8556 & 0.8399 \\
         $G_{23}$ (GPa)  & 0.4494 & 0.4496 \\
         $G_{13}$ (GPa) & 0.8618 & 0.8641 \\
         $\nu_{12}$ & 0.3335  & 0.3337  \\
         $\rho$ (kg/m$^3$) & 1038 & 1036.6 \\
         $\theta_1 (\degree)$ & 0.543 & 0.258 \\
         $\theta_2 (\degree)$ & 44.928 & 44.525 \\
         $\theta_3 (\degree)$ & -45.075 & -45.480 \\
         $\theta_4 (\degree)$ & 90.543 & 90.258 \\ 
         \hline\noalign{\smallskip}
    \end{tabular}
    \label{tab:MPP_validation_Tawfik}
\end{table}

\begin{table}[htbp!]
\captionsetup{justification=centering}
    \caption{Comparison of the safety index $(\beta)$ and probability of failure $(P_f)$ for SSSS $[0\degree/45\degree/-45\degree/90\degree]$ square laminate.}
    \centering
    \begin{tabular}{c c c c}
    \hline\noalign{\smallskip}
         Method &    &  $\beta$ & $P_f$ \\ 
         \hline\noalign{\smallskip}
         \multirow{2}{*}{FORM} &  Present Work & 1.3545 & 0.0878 \\
                &  Tawfik.et.al \cite{tawfik2018neural} &  1.3248 & 0.0926 \\  
         \hline\noalign{\smallskip}
         \multirow{2}{*}{SORM} &  Present Work  & 1.3519 & 0.0882 \\
                &  Tawfik.et.al \cite{tawfik2018neural} &  1.4187 & 0.0780 \\  
         \hline\noalign{\smallskip}
    \end{tabular}
    \label{tab:Pf_beta_comparison}
\end{table}

\subsection{Probability of failure and and global reliability sensitivity analysis for VSCL plate with cutouts}

Next, the Neural Network based probability of failure prediction using the FORM/SORM is extended to VSCL plates with cutouts. As a note to the reader, no benchmark reliability studies of Variable Stiffness Composite Laminate plates based on its fundamental frequency, existed in the literature at the time of writing this paper. To begin with, table \ref{tab:20} summarizes the input parameters that are treated to be random for the case of a 3 ply VSCL composite plate (a/b = 1) along with their mean values, Coefficient of Variation (CoV) and distribution type. As it can be seen, 4 random variables associated with material, 4 w.r.t the cutout shape, size and location and 3 w.r.t the ply angle and its thickness are considered for the reliability analysis. 

A lognormal distribution is considered for all the material parameters and the ply thickness as they it is not logical for them to be assigned with a negative value. Moreover, tighter control in dispersion, i.e. lower value of CoV, is assigned for geometric parameters as it would be expected in any manufacturing setting. Since the FORM-SORM formulation is based on non-correlated standard normal space, any lognormal random variables need to be first converted to its normal space using the following relation: 

\begin{align*}
    \mu_{X} &= 2 ln(\mu_{Y}) - (1/2) ln(\sigma_Y^{2} + \mu_Y^{2})   \\ 
    \sigma_{X}^2 &= -2 ln(\mu_Y) + ln ( \sigma_Y^2 + \mu_Y^2) 
\end{align*}

\noindent where subscript $X$ and $Y$ refers to the random variables in the normal and lognormal space respectively, while $\mu_{(\cdot)}$ and $\sigma_{(\cdot)}$ stands for the mean and standard deviation. 

\begin{table}[htpb!]
\centering
\begin{threeparttable}
\captionsetup{justification=centering}
\caption{Random variables in the analysis of the limit state function w.r.t the first mode of vibration for VSCL plate with cutout.}
\label{tab:20}       
\begin{tabular}{c l c c c c} 
\hline\noalign{\smallskip}
 \multicolumn{6}{c}{Material uncertainties} \\
\hline\noalign{\smallskip}
Sl.No. & Random variables & Symbols & Mean Value & CoV & Distribution type \\
\noalign{\smallskip}\hline\noalign{\smallskip}
1. & Longitudinal modulus (Pa) & $E_{11}$ & 1.73e11 & 0.03701 \tnote{2} & LogNormal \\
2. & Transverse modulus (Pa) & $E_{22}$ & 7.2e9 & 0.03571 \tnote{2} & LogNormal\\
3. & In-plane shear modulus (Pa) & $G_{12}$ & 3.76e9 & 0.05977 \tnote{2} & LogNormal \\
4. & Density (kg/m$^3$) & $\rho$ & 1540 & 0.036 \tnote{2} & LogNormal \\ 
\noalign{\smallskip}\hline
\multicolumn{6}{c}{Geometric uncertainties} \\
\hline\noalign{\smallskip}
5.  & Cutout major axis (m) & $d$ & 0.4 & 0.00025 \tnote{3} & Normal \\
6.  & Cutout ellipticity (-) & $c/d$ & 1 & 0.005 \tnote{3} & Normal\\
7.  & Cutout center location (m) & $x_c$ & 0.0 & 0.001 \tnote{1} & Normal \\
8.  & Cutout center location (m) & $y_c$ & 0.0 & 0.001 \tnote{1} & Normal \\ 
9.  & Ply thickness (m)  & $t^k$ & 0.0033  & 0.04 \tnote{2} & LogNormal \\
10. & Ply angle (\degree) & $<\Theta_0^k | \Theta_1^k>$ & \begin{tabular}{@{}c@{}} $<0\degree,45\degree>$ \\ $<-45 \degree,-60\degree>$ \\ $<0\degree,45\degree>$ \end{tabular} & 1.8\degree \tnote{1} & Normal \\
\hline\noalign{\smallskip}
\end{tabular}
\begin{tablenotes}
\item[1] Standard deviation values are defined for zero mean valued parameters. 
\item[2] CoV used in \cite{tawfik2018neural}.
\item[3] Internal communication with an expert in composite machining. 
\end{tablenotes}
\end{threeparttable}
\end{table}

For the user defined set frequency in Equation \ref{eq:12}, a value of $\lambda_r = 0.97 \lambda_{mean}$ is chosen, where $\lambda_{mean}$ is the deterministic value of first frequency mode of the VSCL plate computed using the mean value for all parameters and plate geometries. 

\subsubsection{ANN based reliability analysis by Monte Carlo simulation} \label{ANN_MCS}

\begin{figure}[htbp!]  
\centering
\includegraphics[width=0.7\textwidth]{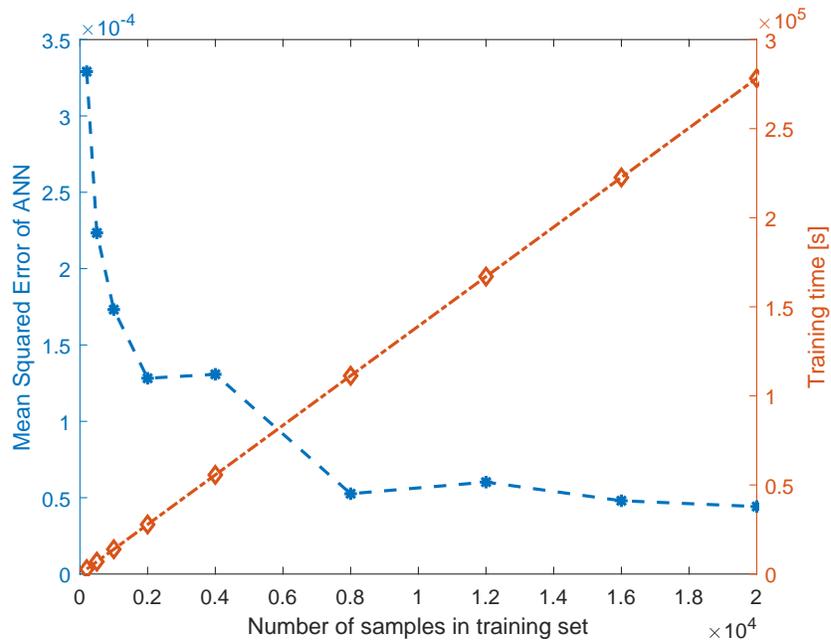}
\captionsetup{justification=centering}
\caption{Time spent on the training set preparation and ANN accuracy as function of the size of the training set.}
\label{fig:ANN_MCS_comp}
\end{figure}

One of the first approaches in using the ANN for reliability estimation is in approximating the Limit State Function (LSF), followed by MCS for reliability assessment. We consider data sets of different sample sizes, ranging from 2000 to 20000 samples, so as to check their influence on the capacity of ANN in approximating the LSF and the reliability analysis error. The training samples are generated using Latin Hypercube Sampling (LHS). \fref{fig:ANN_MCS_comp} plots both the mean square error in ANN approximation of limit state function and the total computational time involved in generating the samples and training the ANN metamodel against the sample size. 

\begin{figure}[htbp!]  
\centering
\includegraphics[width=0.7\textwidth]{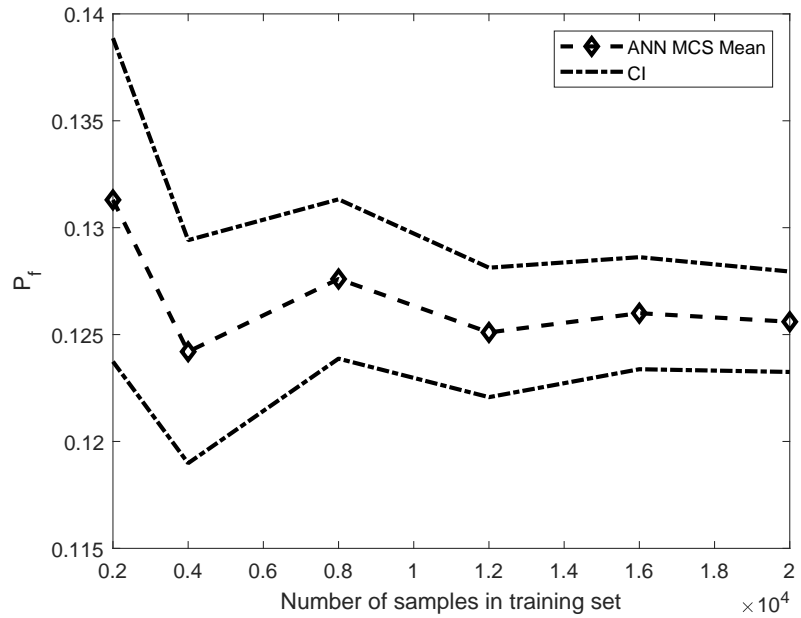}
\captionsetup{justification=centering}
\captionsetup{justification=centering}
\caption{Estimated probability of failure as a function of the number of elements in the training data set. CI stands for confidence interval.}
\label{fig:ANN_MCS_Pf}
\end{figure}

\noindent It can be observed that the error in the network approximation for LSF approaches to zero with increasing number of samples in the input data set. However, larger the data set, more computational time is required for generating those samples. To perform the ANN based reliability studies using MCS, we choose the ANN trained using $8000$ samples, beyond which no significant improvement in ANN approximation of LSF was observed. \fref{fig:ANN_MCS_Pf} plots the estimated failure probability along with the Confidence Interval (CI is defined as twice the standard deviation) (c.f. Equation\ref{eq:1900}) as a function of number of sampling points. As one can observe, the CI turns out to be narrower as the number of samples in the training set increases, due to the effect of reducing standard error of the probability of failure estimations by MCS. From \fref{fig:ANN_MCS_Pf}, it can be observed that a converged failure probability prediction of 0.1251 is obtained using an ANN trained with a data set of 12000 samples, within a prediction uncertainty of 2.4$\%$. 

\subsubsection{ANN based reliability analysis by Monte Carlo based Importance Sampling}

\begin{figure}[htbp!]  
\centering
\includegraphics[width=1\textwidth]{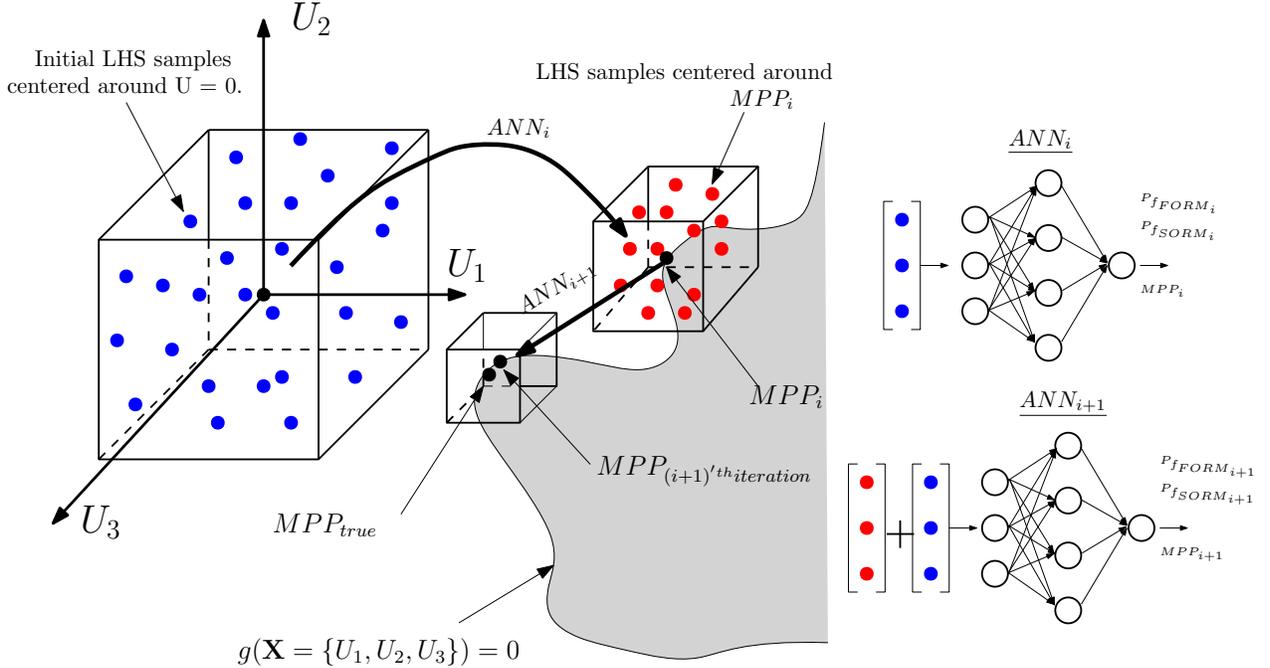}
\captionsetup{justification=centering}
\caption{Schematic sketch showing the flow diagram of the proposed Adaptive ANN based MCIS.}
\label{fig:Adaptive_ANN}
\end{figure}

One of the main drawbacks of ANN based MCS is that it's efficiency depends on the magnitude of the failure probability, meaning, smaller the failure probability, larger the number of samples required. However, this drawback of MCS based reliability can be overcome by forcing the sampling of experimental points using a different probability density function that is centered around the failure domain. This is the fundamental concept behind importance sampling, which was briefly introduced in Section \ref{IS}. Here, we plan to implement an ANN for failure probability estimation that gets trained adaptively using samples chosen from an importance sampling density centered around the MPP estimates in each iterative iterative FORM/SORM step. To make the steps involved more clearer to the reader, as a first step, the network with a given architecture gets trained with a much reduced sample size generated using Latin Hypercube Sampling, centered around the mean values of input parameters. We choose three different training sample set size, namely $N_{samples} = 50,100 ,200$, to train ANN in each iterations. The trained ANN in first stage is then used to find the MPP in the first iteration using FORM/SORM as detailed in Section \ref{ANN_VSCL_prob}. In the second step, an additional $N_{samples}$ samples in the surroundings of the previous MPP from the first iteration are generated using Latin Hypercube Sampling and added to the training sample set and the procedure gets repeated. The process gets terminated once the MPP generated at the i'th iteration within a hypercube that was centered at the MPP from the (i-1)'th iteration and with a width that equals the standard deviation of input random parameters. Once the ANN has been iteratively trained, failure probability is finally computed using the above trained ANN coupled with MCIS around the most recent MPP. Figure \ref{fig:Adaptive_ANN} plots a schematic sketch of the flow diagram of the proposed adaptive ANN based MCIS. Table \ref{tab:ANN_MCIS} presents the results obtained using the above approach. 

\begin{table}[htbp!]
    \captionsetup{justification=centering}
    \caption{Results of ANN based MCIS}
    \centering
    \begin{tabular}{c c c c c c c}
    \hline\noalign{\smallskip}
         \multirow{2}{*}{\makecell{Training set \\ (Step1+Step2)}} & \multicolumn{2}{c} {ANN-training} & \multicolumn{2}{c} {\makecell{$P_{f}$\\(FORM/SORM+MCIS)}} & \multirow{2}{*}{\makecell{Total \\ time [s]}} & \multirow{2}{*}{ANN MSE} \\
         \cline{2-5}  
          & $P_{f_{FORM}}$ & $P_{f_{SORM}}$ & $P_{f_{MCIS}}$ & $Std dev_{P_{f_{MCIS}}}$ &  \\
          \hline 
          50+50 & 0.0924 & 0.0912 & 0.1012 & 0.0016 & 1290 & 3.47$\times 10^{-4}$ \\
          100+100 & 0.1695 & 0.1698 & 0.1152 & 0.0013 & 2587 & 1.67$\times 10^{-4}$ \\
          200+200 & 0.1753 & 0.1736 & 0.1156 & 0.0013 & 5076 & 1.25$\times 10^{-4}$ \\
          \hline\noalign{\smallskip}
    \end{tabular}
    \label{tab:ANN_MCIS}
\end{table}

In the current case, the design point was found to converge within the first two steps of iteration. The estimated failure probability converges to a value of $P_{f_{MCIS}}= 0.1156$ with a variance of $(0.0013)^2$ computed using \eref{eq:1555}. By adopting an adaptive ANN+MCIS approach, we see that the total computational time has been reduced to 5076 seconds, a reduction factor of roughly $95\%$ when compared against the procedure adopted in Section \ref{ANN_MCS}. 

\begin{figure}[htbp]  
\centering
\includegraphics[width=0.9\textwidth]{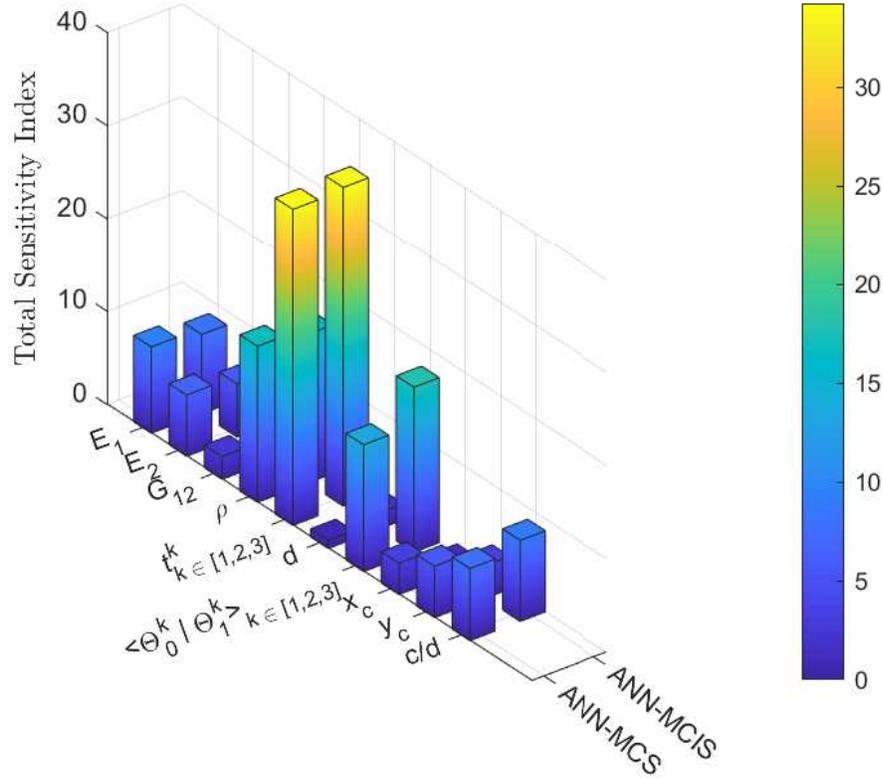}
\captionsetup{justification=centering}
\caption{Estimated total sensitivity index for each of input parameters on the failure probability of aVSCL composite plate. ANN-MCS is based on the ANN metamodel trained using 8000 samples, while ANN-MCIS is based on the metamodel obtained using the proposed adaptive importance sampling method.}
\label{fig:SI_comparison}
\end{figure} 

\fref{fig:SI_comparison} plots the total sensitivity index of the failure probability to the input random parameters. In the figure, the sensitivity indices for each of the ply thickness and for the ply fibre angle orientation both at the mid and edge locations have been coupled together, respectively. The results corresponding to ANN-MCS is obtained using the ANN metamodel trained using 8000 samples (see Section \ref{ANN_MCS}), while the one corresponding to ANN-MCIS are obtained using the ANN metamodel that was adaptively trained, as detailed above. One can observe that an excellent agreement is observed between both the predicted total sensitivity indices for each of the input random parameters. Moreover, the failure probability was found to be the most sensitive on upon the randomness in input ply thickness, followed by ply angle and density uncertainties. 
\subsubsection{Application of Adaptive ANN based MCIS for failure prediction and reliability sensitivity studies for different composite plates} \label{ANN_comparison}

In this section, we extend our proposed adaptive ANN based MCIS for failure prediction and reliability estimates of four different composite plates. Table \ref{tab:Diff_composites} summarizes four different types of VSCL composites considered for this study, with variations in aspect ratio and ply angle orientations. 

\begin{table}[htbp!]
    \captionsetup{justification=centering}
    \caption{Different types of VSCL composites identified for failure probability and reliability estimates using proposed adaptive ANN-MCIS.}
    \centering
    \begin{tabular}{c c c c}
    \hline\noalign{\smallskip}
         Composite Type & Aspect ratio (a/b) & Ply angle orientation & Remarks\\
         \hline\noalign{\smallskip}
         Composite 1 & a/b = 1 & \begin{tabular}{@{}c@{}} $<0\degree,45\degree>$ \\ $<-45 \degree,-60\degree>$ \\ $<0\degree,45\degree>$ \end{tabular} & \makecell{Square plate \\ with unsymmetric ply} \\
         \hline\noalign{\smallskip}
         Composite 2 & a/b = 1/2 & \begin{tabular}{@{}c@{}} $<0\degree,45\degree>$ \\ $<-45 \degree,-60\degree>$ \\  $<0\degree,45\degree>$ \end{tabular} & \makecell{Rectangular plate \\ with unsymmetric ply} \\ 
          \hline\noalign{\smallskip}
         Composite 3 & a/b = 1 & \begin{tabular}{@{}c@{}} $<0\degree,45\degree>$ \\ $<-45 \degree,-60\degree>$ \\  $<-45 \degree,-60\degree>$ \\ $<0\degree,45\degree>$ \end{tabular} & \makecell{Square plate \\ with symmetric ply} \\
         \hline\noalign{\smallskip}
         Composite 4 & a/b = 1/2 & \begin{tabular}{@{}c@{}} $<0\degree,45\degree>$ \\ $<-45 \degree,-60\degree>$ \\  $<-45 \degree,-60\degree>$ \\ $<0\degree,45\degree>$ \end{tabular} & \makecell{Rectangle plate \\ with symmetric ply} \\
         \hline\noalign{\smallskip}
    \end{tabular}
    \label{tab:Diff_composites}
\end{table}

Table \ref{tab:ANN_MCIS_4_composites} summarises the failure probabilities of the identified 4 different types of composites w.r.t the given input random parameters (c.f. Table \ref{tab:20}). The failure probabilities were computed using the proposed adaptive ANN based MCIS, which was trained using 200 samples for each iteration. 

\begin{table}[htbp!]
    \captionsetup{justification=centering}
    \caption{Results of failure probabilities of four different types of VSCL plate configurations using the proposed adaptive ANN based MCIS.}
    \centering
    \begin{tabular}{c c c c c c}
    \hline\noalign{\smallskip}
         \multirow{2}{*}{\makecell{Composite \\ Type}} & \multirow{2}{*}{\makecell{$F_{det}$ \\ $[rad/s]$ }} & \multicolumn{2}{c} {ANN-training} & \multicolumn{2}{c} {\makecell{$P_{f}$\\(FORM/SORM+MCIS)}} \\
         \cline{3-6}  
          & & $P_{f_{FORM}}$ & $P_{f_{SORM}}$ & $P_{f_{MCIS}}$ & $Std dev_{P_{f_{MCIS}}}$ \\
         \hline 
          1 & 1193.5 & 0.1753 & 0.1736 & 0.1156 & 0.0013 \\
          2 & 406.08 & 0.1898 & 0.1901 & 0.1897 & 0.0016 \\
          3 & 1820.4 & 0.1674 & 0.1690 & 0.1014 & 0.0012 \\
          4 & 589.06 & 0.1842 & 0.1843 & 0.1813 & 0.0015 \\
          \hline\noalign{\smallskip}
    \end{tabular}
    \label{tab:ANN_MCIS_4_composites}
\end{table}

The second column of Table \ref{tab:ANN_MCIS_4_composites} contains information on the deterministic fundamental frequency of respective composite types, computed at mean values of input random parameters. A value which is $97\%$ of these deterministic values serves as $\lambda_r$ in \eref{eq:12}, which defines there respective limit-state function definition. 

\begin{figure}[htbp!]  
\centering
\includegraphics[width=0.9\textwidth]{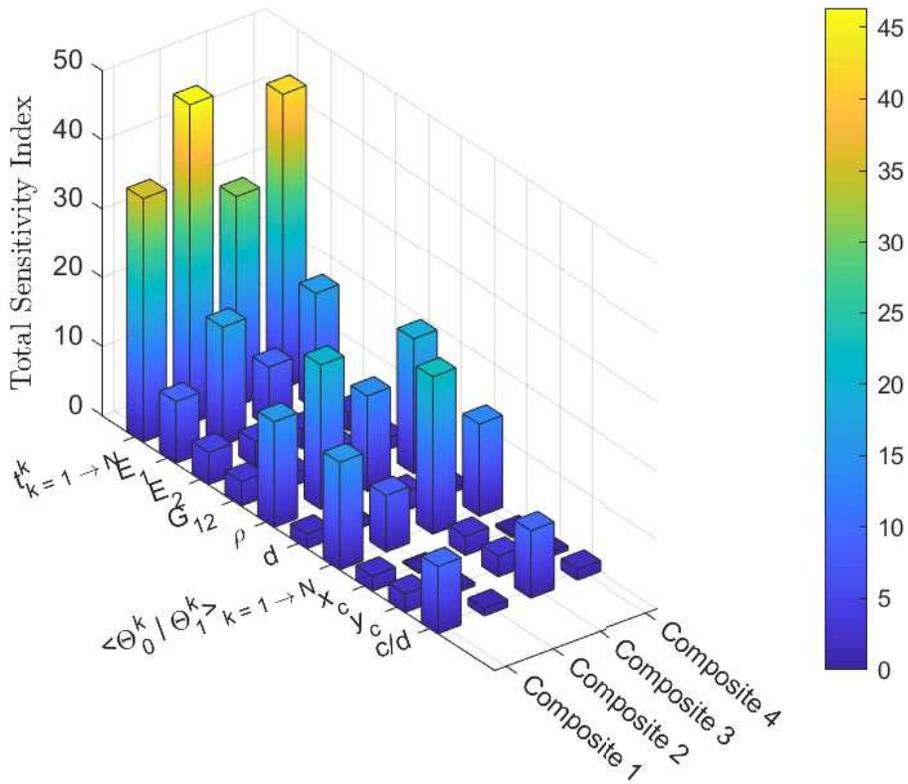}
\captionsetup{justification=centering}
\caption{Estimated total sensitivity index for each of input parameters on the failure probability for four different VSCL composite plate (c.f. Table \ref{tab:Diff_composites}).}
\label{fig:SI_comparison_4_composites}
\end{figure} 

Figure \ref{fig:SI_comparison_4_composites} plots the total sensitivity index of each of the input parameters on the failure probability, for each of the four different VSCL plate configurations considered. The total sensitivity index are computed using the adaptive ANN based MCIS metamodel. From Table \ref{tab:ANN_MCIS_4_composites}, it is clear that for a given input random parametric data, the rectangular plate with unsymmetric ply layup was found to show the maximum probability of failure, followed closely by its symmetric ply layup counterpart. Square plates, irrespective of the ply layup, was found to exhibit a lower failure probability by roughly a factor of $40\%$, compared to its rectangular counterpart. 

Furthermore, by having a closer look at Figure \ref{fig:SI_comparison_4_composites}, we first see that the failure probability remains least sensitive to changes in certain input random parameters, irrespective of the VSCL plate configurations, namely, $E_2, G_{12}, d, x_c $ and $y_c$ (Total sensitivity index below $5\%$). On the other spectrum of sensitivity values, we see that the total uncertainty in ply thickness dominates in its contribution towards the failure probability, across all the 4 VSCL configurations considered. This observation matches very well with the results published in the work of Tawfik\textit{.et.el.} \cite{tawfik2018neural}, where the importance of taking into account the uncertainty of ply thickness on failure estimates of straight fiber composite case studies was emphasised. Moreover, changing the plate configuration from a square to a rectangular configuration was only found to exacerbate the contribution of uncertain ply thickness to the overall failure probability. However, the reverse trend was observed in the failure sensitivity w.r.t the ply angle orientation, such that for a given ply layup (un/symmetric), play angle uncertainties in square plates were found to be more sensitive to the failure probability than its rectangular counterpart. 

Among the material parametric uncertainties, only $E_1$ and $\rho$ was found to have a significant contribution, where the later was found to consistently exhibit greater impact on the failure probability over the later across the VSCL configurations. Finally, as far as the effect of uncertainty in hole parameters on the failure probability, its ellipticity ($c/d$) was found to be the most dominant of them all, especially in the case of square plates.

\section{Conclusion} \label{Conclusion}

In this paper, both the failure probability estimates and sensitivity studies of various configurations of VSCL plates for a given set of random input parametric values has been thoroughly discussed, analyzed and its performance validated with that of conventional MCS. The geometric uncertainties in the cutout geometry was incorporated into the studies using the stochastic level set method. Aside from the conventional methods of ANN based reliability estimates for composites where randomly large number of samples are drawn from the input parametric space, our proposed method of adaptive ANN based Importance Sampling shows promise in delivering accurate results for failure predictions and total sensitivity indices, but at a cost of significant computational savings by an observed factor of nearly $95\%$ compared to that of MCS. Moreover, the importance of having to take into account the uncertainties in ply thickness for failure probability estimates is also highlighted, irrespective of the plate aspect ratio or the ply layup sequence. 


\section*{References}
\bibliographystyle{model1-num-names}
\bibliography{Reference}
\end{document}